\documentclass[MSNbibl,number,citesort,dvips]{arxbj}
\usepackage{mathbh}
\usepackage{dcolumn}
\usepackage{mathrsfs}
\usepackage{graphicx}

\aid{0}
\volume{19}
\issue{3}
\pubyear{2013}
\firstpage{905}
\lastpage{930}
\doi{10.3150/12-BEJ419} 

\makeatletter
\newcolumntype{d}[1]{D{.}{.}{#1}}
\newproclaim{definition}{Definition}
\newtheorem{proposition}{Proposition}
\newremark{remarque}{Remark}

\newcommand{\eqref}[1]{(\ref{#1})}
\newcommand{\EE}{\mathbf{E}}
\newcommand{\redpalm}[1]{P^{!{#1}}}
\newcommand{\Eredpalm}[1]{\EE^{!{#1}}}
\newcommand{\R}{\mathbb{R}}
\newcommand{\Rd}{\mathbb{R}^d}
\newcommand{\Z}{\mathbb{Z}}
\newcommand{\Zd}{\mathbb{Z}^d}
\newcommand{\G}{\mathcal{G}}
\newcommand{\Hfamily}{\mathcal H}
\newcommand{\HH}{\mathbf{H}}
\newcommand{\Oi}{\Omega_{\neq\infty}}
\newcommand{\dee}{\,\mathrm{d}}
\newcommand{\1}{\mathbh{1}}
\renewcommand{\div}{\operatorname{div}}

\newcommand{\xx}{\mathbf x}
\newcommand{\yy}{\mathbf y}
\newcommand{\zz}{\mathbf z}
\newcommand{\transpose}[1]{{#1}^\top}
\newcommand{\coord}[2]{{#1}^{(#2)}}
\newcommand{\listelem}[2]{{\mathbf{#1}}_{#2}}
\newcommand{\coordlistelem}[3]{{#1}_{#2}^{(#3)}}
\newcommand{\listof}[2]{{\mathbf{#1}}_{#2}}
\newcommand{\bs}{\setminus}
\newcommand{\C}{\mathcal{C}}
\newcommand{\D}{\Delta}
\renewcommand{\L}{\Lambda}
\newcommand{\LL}{\mathcal{R}}
\newcommand{\om}{\omega}
\renewcommand{\O}{\Omega}
\renewcommand{\t}{\theta}
\renewcommand{\tt}{{\boldsymbol{\theta}}}
\newcommand{\V}{\mathcal{V}}
\newcommand{\vp}{\varphi}
\newcommand{\vs}{\varsigma}
\makeatother

\begin{document}
\begin{frontmatter}

\title{Variational estimators for the parameters of Gibbs point
process models}
\runtitle{Variational estimators}

\begin{aug}
\author[1,2]{\fnms{Adrian} \snm{Baddeley}\corref{}\thanksref{1,2}\ead[label=e1]{Adrian.Baddeley@csiro.au}} \and
\author[3]{\fnms{David} \snm{Dereudre}\thanksref{3}\ead[label=e2]{david.dereudre@univ-lille1.fr}}
\runauthor{A. Baddeley and D. Dereudre} 
\address[1]{CSIRO Mathematics, Informatics and Statistics,
Leeuwin Centre, 65 Brockway Rd, Floreat WA 6014,
Australia. \printead{e1}}
\address[2]{School of Mathematics \& Statistics, University of Western Australia}
\address[3]{Laboratoire Paul Painlev\'e, UMR CNRS 8524, University of Lille, 59655 Villeneuve d'Ascq C\'edex,
France. \printead{e2}}
\end{aug}

\received{\smonth{8} \syear{2010}}
\revised{\smonth{9} \syear{2011}}

%
\begin{abstract}
This paper proposes a new estimation technique for fitting parametric
Gibbs point process models to a spatial point pattern dataset.
The technique is a counterpart, for spatial point processes,
of the variational estimators for Markov random fields
developed by Almeida and Gidas.
The estimator does not require the point process density
to be hereditary, so it is applicable to models which do not
have a conditional intensity, including models which exhibit
geometric regularity or rigidity.
The disadvantage is that the intensity
parameter cannot be estimated: inference is effectively
conditional on the observed number of points.
The new procedure is faster and more stable than
existing techniques, since it does not require simulation, numerical
integration or optimization with respect to the parameters.
\end{abstract}

%

\end{frontmatter}

\section{Introduction}
Statistical methodology for fitting models to spatial point pattern data
has been a subject of intensive research for three decades.
Likelihood-based methods were
once regarded as computationally prohibitive (Ripley~\cite{ripl88}, Introduction).
Although maximum likelihood and Bayesian methods can now be implemented
using Markov chain Monte Carlo (Geyer~\cite{geye99}), this approach is still
computationally intensive. Various alternative strategies have been
explored, including analytic approximations
to the likelihood (Ogata and Tanemura~\cite{ogattane84}),
computationally efficient surrogates such as maximum pseudolikelihood
(composite likelihood) (Besag~\cite{besa78})
and Takacs--Fiksel estimators (Fiksel~\cite{fiks84}, Takacs~\cite{taka86}).

These approaches fail when the point pattern data exhibit a high
degree of geometric regularity or rigidity. For example, point patterns
which approach a random dense packing of hard spheres, or
a randomly perturbed hexagonal lattice, can be constructed as realisations
of certain Gibbs models. Existing methods for fitting Gibbs models
to the data, such as maximum pseudolikelihood
(Baddeley and Turner~\cite{baddturn00}, Besag~\cite{besa78}, Billiot \textit{et~al.}~\cite{billcoeudrou08}, Coeurjolly and Drouilhet~\cite{coeudrou10}, Goulard \textit{et~al.}~\cite{goulsarkgrab96}),
tend to be numerically unstable in the nearly-rigid case.
The likelihood-based procedures require very long computation time
for simulations. Other existing methods are generally based on
an equilibrium equation involving the addition or removal of points
of the process; such transitions may be impossible, or rare, if the model
is too rigid. Moreover, most of the existing methods are based on the
(Papangelou) conditional intensity; in some recent work (Dereudre~\cite{Dereudre},
 Dereudre and Lavancier~\cite{DL}, Dereudre \textit{et~al.}~\cite{DDG}),
geometrically rigid point patterns are generated using Gibbs models which
violate the usual assumption that the probability density is hereditary,
so that the conditional intensity may not exist.
Although the classical procedures of pseudolikelihood and
Takacs--Fiksel estimation have been generalized to the non-hereditary
setting (Coeurjolly \textit{et~al.}~\cite{CDDL}, Dereudre and Lavancier~\cite{DL}), the associated estimators remain unavailable
or inefficient if the process is too rigid (see, e.g., the simulations
presented in Dereudre and Lavancier~\cite{DL2}).

In this paper, we propose an alternative approach to parameter estimation
which is motivated by the variational estimators of Almeida and Gidas
\cite{AG} for discrete space Markov random fields.
This approach does not require the hereditary property.
It is based on an equilibrium equation involving infinitesimal
perturbation of the local energy and test functionals (see equation
\eqref{equationvs}). The use of infinitesimal perturbations seems
naturally well-adapted for rigid models, since it does not require
addition or removal of points. A necessary assumption to obtain
identification is that the interaction potential is not
constant. Therefore, one of the first consequences is that the intensity
parameter $z$ of the Gibbs process cannot be fitted, since it
corresponds to a constant point potential $-\ln(z)$.
Let us note that the parameters of a Gibbs
process (apart from $z$) can be estimated without knowing the value of
$z$. If necessary, the intensity parameter $z$ may be fitted in a
second step
using another procedure, such as maximum pseudolikelihood, which performs
well when interaction parameters are fixed.

Finally, we note that our procedure is quicker than
existing ones, since it does not require simulation, numerical integration
or optimization with respect to the parameters. The
algorithm is very simple to implement, and requires only the
computation of sums and the inversion of a linear system of size equal
to the
number of parameters (see~(\ref{estimator}),~(\ref{estimatorb})).
The estimator is exact and explicit.

In Section~\ref{S:defns}, we introduce notation and basic definitions
for (grand canonical) Gibbs processes.
In Section~\ref{SectionVE}, the
variational equilibrium equations are stated in stationary and non-stationary
versions. In Section~\ref{S:expo}, the variational procedure is also
presented in two versions, corresponding to the two variational
equations. In Section~\ref{SectionA}, asymptotic properties are
investigated. We show that both procedures are strongly consistent
and we prove asymptotic normality for one of them.
In Section~\ref{SectionE},
we present a large class of examples for which the procedure is
available. We choose two typical examples coming from
statistical mechanics and stochastic geometry (Lennard-Jones model,
Hard sphere model). Section~\ref{SectionS}
is devoted to simulation experiments, where the variational procedure
is applied
in three situations corresponding to non-rigid, rigid and very rigid
cases of the Lennard-Jones point process.
%
\section{Definitions and notation}\label{S:defns}\vspace*{-1pt}
%
\subsection{State spaces and reference measures}\vspace*{-1pt}
Our setting is Euclidean space $\Rd$ of arbitrary dimension $d \ge1$.
An element of $\Rd$ is denoted by
$\xx= (\coord x 1, \ldots, \coord x d )$.
Lebesgue measure on $\Rd$ is denoted by ${\lambda}^d$.
\begin{definition}
A \textup{configuration} is a subset $\om$ of $\Rd$
which is locally finite, meaning that $\om\cap\L$ has finite cardinality
$N_\L(\om)=\#(\om\cap\L)$ for every bounded Borel set $\L$.
The space $\O$ of all configurations is equipped with the
$\sigma$-algebra $\mathcal{F}$ generated by the counting variables~$N_\L$.
\end{definition}

The symbol $\L$ will always refer to a bounded Borel set in $\R^d$.
It will often be convenient to write $\om_\L$ in place of $\om\cap
\L$.
We abbreviate $\om\cup\{\xx\}$ to $\om\cup\xx$
and abbreviate $\om\backslash\{\xx\}$ to $\om\backslash\xx$
for every $\om$ and every $\xx$ in $\om$.
For $k$ points $\xx_1,\xx_2,\ldots,\xx_k$ in $\Rd$,
we denote by $\xx_{1\ldots k}$ the configuration $\xx_1\cup\cdots
\cup\xx_k$
and for $1\le i \le k$, denote by $\xx_{1\ldots k\bs i}$ the
configuration $\xx_{1\ldots k} \bs\xx_i$.

As usual, we take the reference measure on $(\O,\mathcal{F})$
to be the distribution $\pi^z$ of the Poisson point process
with intensity measure $z{\lambda}^d$ ($z>0$) on $\Rd$.
Recall that $\pi^z$ is the unique probability measure
on $(\O,\mathcal{F})$ such that the following hold for all subsets
$\L$:
(i) $N_{\L}$ is Poisson distributed with parameter $z
{\lambda}^d(\L)$,
and (ii) conditional on $N_\L=n$,
the $n$ points in $\L$ are independent with uniform distribution on
$\L$,
for each integer $n \ge1$.
The Poisson point process restricted to $\L$ will be denoted $\pi
^z_\L$.

In this paper, we consider only point process distributions
(probability measures $P$ on~$\O$)
such that the intensity measure $m$ on $\Rd$, defined by $m(\L)=\EE
_P(N_\L)$,
for any~$\L$, is $\sigma$-finite.
Here $\EE_P$ denotes expectation with respect to $P$.
We denote by $C_P^!$ the reduced Campbell measure of $P$
defined on $\Rd\times\O$ by
%
\begin{equation}
C_P^!(g)=\int\sum_{\xx\in\om} g(\xx,\om\bs\xx) P(\mathrm{d}\om)
\end{equation}
for any positive measurable function $g$ from $\Rd\times\O$ to $\R$.

Translation by a vector $u \in\Rd$ is denoted by $\tau_u$,
whether acting on $\Rd$ or on $\O$.
When $P$ is stationary (i.e., $P=P\circ\tau_u^{-1}$ for any $u$ in
$\Rd$)
the intensity measure has the form $m = z(P){\lambda}^d$ and
there exists a unique probability measure $\redpalm0$,
called the reduced Palm measure, such that
%
\begin{equation}
\label{Palm}
C_P^!(g)=z(P)\int\int g(\xx,\tau_\xx\om) \redpalm0(\mathrm{d}\om
){\lambda}^d(\mathrm{d}\xx).
\end{equation}
See Matthes \textit{et~al.}~\cite{KMM} for more details about Campbell and Palm measures.\vspace*{-1pt}
%
\subsection{Interaction}\vspace*{-1pt}
We shall define the interaction energy in a general setting,
along the lines of Preston~\cite{Preston}.
Thus, we do \emph{not} assume that the local densities come from a multibody\vadjust{\goodbreak}
interaction potential.
This general viewpoint allows us to deal with the non-hereditary case
(see Definition~\ref{hereditary} below).
\begin{definition}
A \textup{family of energies} is a collection $\Hfamily= (H_\L)$,
indexed by bounded Borel sets $\L$,
of measurable functions from $\O$ to $\R\cup\{+\infty\}$ such that,
for every $\L\subset\L'$,
there exists a measurable function
$\varphi_{\L,\L'}$ from $\O$ to $\R\cup\{+\infty\}$
such that for every $\om\in\O$
%
\begin{equation}\label{compatible}
H_{\L'}(\om)= H_{\L}(\om) + \varphi_{\L,\L'}(\om_{\L^c}).
\end{equation}
\end{definition}

Equation (\ref{compatible})
is equivalent to (6.11) and (6.12) in Preston~\cite{Preston}, page 92.
In physical terms, $H_\L(\om)=H_\L(\om_\L\cup\om_{\L^c})$ represents
the potential energy of $\om_\L$ inside $\L$ given
the configuration $\om_{\L^c}$ outside $\L$.
\begin{definition}
A configuration $\om$ has \textup{locally finite energy} with respect to
a family of energies $\Hfamily= (H_\L)$ if, for every $\L$,
the energy $H_{\L}(\om)$ is finite.
We denote by $\Oi(\Hfamily)$ or simply $\Oi$
the space of configurations which have a locally finite energy.
\end{definition}
\begin{definition}
\label{hereditary}
A family of energies $\Hfamily= (H_\L)$ is \textup{hereditary}
if for all $\L$, all $\om\in\O$ and $\xx\in\L$
%
\begin{equation}
\label{her}
H_\L(\om)=+\infty \quad\Rightarrow\quad H_{\L}(\om\cup\xx)=+\infty.
\end{equation}
\end{definition}

The assumption \eqref{her} is necessary in many papers,
for example, Nguyen and Zessin~\cite{NZ}, Ruelle~\cite{Ruelle70}.
In this setting, the local energy $h(\xx,\om)$ is defined
for every $\om\in\Oi$ and $\xx\notin\om$ by
%
\begin{equation}\label{localenergy}
h(\xx,\om)= H_\L(\om\cup\xx)-H_\L(\om)
\end{equation}
for any $\L$ containing $\xx$.
Note that by (\ref{compatible}), this definition does not depend on
$\L$.

Some recent work deals with non-hereditary Gibbs models (Dereudre~\cite{Dereudre}, Dereudre and Lavancier~\cite{DL}, Dereudre \textit{et~al.}~\cite{DDG}).
This setting occurs only if the family of energies has a hardcore part,
that is, if $\Oi((H_\L)) \neq\O$.
Henceforth we will \emph{not} assume that the energy is hereditary.
The equations and estimators presented in the following sections
are available in the hereditary or non-hereditary setting.
%
\subsection{Gibbs point processes}
We are now in a position to define Gibbs measures.
Let us make an integrability assumption on the family of energies,
equivalent to (6.8) in Preston~\cite{Preston}.
\begin{definition}
\label{integrable}
The family of energies $\Hfamily= (H_\L)$ is \textup{integrable} if,
for every $\L$ and every $\om$ in $\Oi$, we have
%
\begin{equation}
\label{integ}
0< \int \mathrm{e}^{-H_\L(\om'_\L\cup\om_{\L^c})}\pi^z_\L(\mathrm{d}\om'_\L
) < +\infty.
\end{equation}
\end{definition}

The second inequality in (\ref{integ}) is in general ensured by the
stability of the energy functions. The first inequality is obvious in
the classical hereditary setting, while in the non-hereditary case, it
remains true under reasonable assumptions (e.g., Dereudre~\cite{Dereudre}, Dereudre \textit{et~al.}~\cite{DDG}).

Under this integrability assumption, for every $\L$ and every $\om$
in $\Oi$,
the local conditional density $f_\L$ is defined by
%
\begin{eqnarray}
\label{localdensity}
f_\L(\om)  =  \frac{1}{Z_\L(\om_{\L^c})} \mathrm{e}^{-H_\L(\om)},
\end{eqnarray}
where $ Z_\L(\om_{\L^c})$ is the normalization constant defined by
$Z_\L(\om_{\L^c})= \int \mathrm{e}^{-H_\L(\om'_\L\cup\om_{\L^c})}\pi
^z_\L(\mathrm{d}\om'_\L)$.
Note that from (\ref{integ}), $ 0< Z_\L(\om_{\L^c})< +\infty$ and
therefore this local density is well-defined.

The usual definition of a ``Gibbs point process'' is
equivalent to the following (Georgii~\cite{geor88}, page 28).
\begin{definition}
A probability measure $P$ on $\O$ is a \textup{(grand canonical)
Gibbs measure}
for the integrable family of energies $\Hfamily= (H_\L)$
and the intensity $z>0$
if $P(\Oi)=1$ and, for every $\L$,
for any measurable and integrable function $g$ from $\O$ to $\R$,
%
\begin{equation}\label{DLR}
\int g(\om)P(\mathrm{d}\om) = \int\int g(\om'_\L\cup\om_{\L^c}) f_\L
(\om'_\L\cup\om_{\L^c}) \pi^z_\L(\mathrm{d}\om'_\L) P(\mathrm{d}
\om).
\end{equation}
Equivalently,
for $P$-almost every $\om$ the conditional law of $P$ given $\om_{\L^c}$
is absolutely continuous with respect to $\pi^z_\L$ with the density
$f_\L$.
\end{definition}

The equations (\ref{DLR}) are called the Dobrushin--Lanford--Ruelle (DLR)
equations.
We denote by $\G$ the set of Gibbs measures.

For every $\L$, every $k\ge0$ and every $\om$ in $\Oi$,
define
\[
Z_{\L,k}(\om_{\L^c}) =
\int_{\L^k} \mathrm{e}^{-H_\L(\xx_{1 \ldots k} \cup\om_{\L^c})}
\dee\xx_1\cdots\dee\xx_k.
\]
Under the integrability assumption \eqref{integ},
$Z_{\L,k}(\om_{\L^c})$ is always finite and there exists
at least one $k$ such that $Z_{\L,k}(\om_{\L^c})>0$.
Provided $0<Z_{\L,k}(\om_{\L^c})<\infty$ we may define
the local conditional density for fixed $k$,
%
\begin{equation}
\label{localdensity.grand}
f_{\L,k}(\om)  =  \frac{1}{Z_{\L,k}(\om_{\L^c})} \mathrm{e}^{-H_\L(\om)}
\1\{N_\L(\om) = k\}.
\end{equation}

\begin{definition}\label{GCG}
A probability measure $P$ on $\O$ is a (canonical) Gibbs measure
for the integrable family of energies $\Hfamily= (H_\L)$ if $P(\Oi
)=1$ and,
for every $\L$,
for any measurable and integrable function $g$ from $\O$ to $\R$,
%
\begin{eqnarray}\label{DLRgc}
\int g(\om)P(\mathrm{d}\om) &=&
\sum_{k=0}^\infty\int \int_{\L^k}
\1_{\{N_\L(\om)=k\}}g(\xx_{1\ldots k} \cup\om_{\L^c})\nonumber
\\[-8pt]
\\[-8pt]
&&\hphantom{\sum_{k=0}^\infty\int \int_{\L^k}}\times
f_{\L,k}(\xx_{1\ldots k} \cup\om_{\L^c})
\dee\xx_1\cdots\dee\xx_k P(\mathrm{d}\om).\nonumber
\end{eqnarray}
Equivalently, for every $k\ge0$,
for $P$-almost every $\om$ such that $f_{\L,k}(\om)$ is well defined,
the conditional law of $P$ given $\om_{\L^c}$ and $N_\L(\om)=k$
is absolutely continuous with respect to $\pi_\L(\cdot|N_\L=k)$
with density $|\L|^k f_{\L,k}$.
\end{definition}

The results in Section~\ref{SectionVE} below
are proved for canonical Gibbs measures.
It is obvious that any Gibbs measure $P$ is also a canonical Gibbs
measure. Therefore, the results remain true for Gibbs measures.

Let us note that any canonical Gibbs measure
admits the mixture representation $P = \int P_z \nu(\mathrm{d} z)$, where
$\nu$ is a probability measure on $\R^+$ and $P_z$ a Gibbs measure
for intensity~$z$ (see Georgii~\cite{G} and Theorems 2.1 and 2.2 in Preston~\cite{Preston}).

Let us recall the famous Nguyen--Zessin equation (Nguyen and Zessin~\cite{NZ}), which is available
only in the hereditary setting. For any Gibbs measure $P$,
%
\begin{equation}\label{NZess}
C_P^!=z\mathrm{e}^{-h}{\lambda}^d\otimes P.
\end{equation}
A generalization in the non-hereditary setting is investigated in Dereudre and Lavancier~\cite{DL}.
%
\section{Variational equation for Gibbs point processes}\label{SectionVE}
%
\subsection{Definitions and notation}
First, we introduce the concept of a differentiable point in a configuration.
\begin{definition}
Let $\om$ be a configuration
and $\xx$ a point in $\Rd$ such that $\xx\notin\om$.
We say that a family $\Hfamily= (H_\L)$ is \textup{differentiable}
at $(\xx,\om)$ if there exists an open neighborhood $\V$ around $\xx
$ in $\Rd$
and a bounded set $\L$ containing $\V$ such that the function
$\V\to\R\cup\{\infty\}$ defined by
$\yy\mapsto H_\L(\om\cup\yy)$
is differentiable at $\xx$ in the usual sense. By convention,
if this function is equal to infinity on $\V$, we take the derivative
to be zero. We denote by $\nabla H(\xx,\om)$ the corresponding gradient,
and by $\frac{\partial H}{\partial\coord x i}(\xx,\om)$ its $i$th
coordinate.
\end{definition}

The notation $\nabla H(\xx,\om)$ is well-defined
since this quantity does not depend on $\L$ by~\eqref{compatible}.
Note that, in the hereditary setting, differentiability of $(H_\L)$
is equivalent to differentiability of the local energy $h(\xx,\om)$
with respect to the first variable $\xx$.
Then it is clear that $\nabla H(\xx,\om)=\nabla_\xx h(\xx,\om)$.

We say that $(H_\L)$ is \emph{${\lambda}^d$-a.e.
differentiable} if, for every $\om$,
$\nabla H(\xx,\om)$ exists for ${\lambda}^d$-almost
every~$\xx$.
Henceforth we assume that $(H_\L)$ is ${\lambda}^d$-a.e.
differentiable.

Next, we introduce classes of functionals used in the sequel.
\begin{definition}
Let $g$ be a measurable function from $\Rd\times\O$ to $\R$.

We say $g$ is ${\lambda}^d$\textup{-a.e. differentiable}
if for any $\om$ the
function $\xx\mapsto g(\xx,\om)$ is ${\lambda}^d$-a.e.
differentiable.
We denote the gradient function by $\nabla g(\xx,\om)$ and
its $i$th coordinate by $\frac{\partial g}{\partial\coord x i}(\xx
,\om)$.

We say that $g$ is \textup{shift invariant} if, for any vector $u$ in
$\Rd$,
$g(\tau_u (\xx),\tau_u (\om))=g(\xx,\om)$ for all $\xx$ and $\om$.\vadjust{\goodbreak}

We say that $g$ has \textup{compact support} if there exists a compact set
$K$ in $\Rd$ such that $g(\xx,\om)=0$ for all $\xx$ outside $K$ and all
$\om$.\vspace*{-1pt}
\end{definition}

Recall that a function $f\dvtx \R\to\R$ is absolutely continuous
on $[a,b]$ if it is differentiable ${\lambda}$-a.e. on
$[a,b]$ and
if, for any $x \in[a,b]$, $f(x)-f(a)=\int_a^x f'(t)
{\lambda}(\mathrm{d} t)$.
For example any Lipschitzian function is absolutely continuous.\vspace*{-1pt}
\begin{definition}
Let $g$ be a measurable function from $\Rd\times\O$ to $\R$.
We say $g$ is \textup{regularizing} with respect to $(H_\L)$
if for every $\om$, every $a>0$, every $1\le i \le d$ and $
{\lambda}$-a.e. every
$\coord x 1, \ldots, \coord x {i-1} , \coord x {i+1} , \ldots, \coord
x d$
in $[-a,a]$ the function $[-a,a] \to\R\cup\{\infty\}$ defined by
%
\begin{equation}\label{functl}
\coord x i \mapsto g(\xx,\om)\mathrm{e}^{-H_{[-a,a]^d}(\om\cup\xx)}
\end{equation}
is absolutely continuous.\vspace*{-1pt}
\end{definition}

If $g$ is ${\lambda}^d$-a.e. differentiable and
regularizing, then for ${\lambda}^d$-almost every $\xx$
the derivative of \eqref{functl} is equal to
\[
\biggl(
\frac{\partial g}{\partial\coord x i}(\xx,\om)
-g(\xx,\om)\frac{\partial H}{\partial\coord x i}(\xx,\om)
\biggr)
\mathrm{e}^{-H_{[-a,a]^d}(\om\cup\xx)}.
\]

In many situations, it becomes easier to check the regularizing condition,
in particular when the energy functions $(H_\L)$ are regular.
For instance, in many examples in Section~\ref{SectionE},
the function \eqref{functl} will be continuous and piecewise differentiable
which ensures absolute continuity.
Nevertheless, we state the most general condition
in order to investigate the largest possible class of examples.

We denote by $\C_\LL= \C_\LL(H_\L)$ the class of functions $g$ which
are ${\lambda}^d$-a.e. differentiable and regularizing
with respect to $(H_\L)$.
The subclass $\C_{\LL,K}$ consists of those functions $g \in\C_\LL
$ which
have compact support; the subclass $\C_{\LL,\tau}$ consists of those $g$
which are shift invariant.\vspace*{-1pt}
%
\subsection{The general variational equation}\vspace*{-1pt}
\begin{proposition}\label{propv}
Let $P$ be a canonical Gibbs measure
for the energy functions $(H_\L)$. Then for every function $g$ in the class
$\C_{\LL,K}$ such that
%
\begin{equation}\label{integrability}
C_P^!(|\nabla g|+|g\nabla H|)<\infty,
\end{equation}
we have
%
\begin{equation}\label{equationv}
C_P^!(\nabla g)=C_P^!(g\nabla H).
\end{equation}
\end{proposition}

First, we remark that the expectation in (\ref{integrability}) is meaningful.
Indeed let $g$ be a function in $\C_{\LL,K}$ with compact support $\L$.
From \eqref{DLRgc},
%
\begin{eqnarray}\label{devCamp}
C_P^!(|\nabla g|+|g\nabla H|)& =& \int\sum_{\xx\in\om_\L} |\nabla g(\xx,\om\bs\xx
)|+|g\nabla H(\xx,\om\bs\xx)| P(\mathrm{d}\om) \nonumber\\
&=& \sum_{k=1}^{+\infty} \int\frac{\1_{\{ N_\L(\om)=k \}}}{Z_{\L
,k}(\om_{\L^c})} \sum_{i=1}^k \int_{\L^k} (|\nabla g|+|g\nabla
H|)(\listelem x i, \om_{\L^c}\cup\listof x {1\ldots k\bs i})\\
&&\hphantom{= \sum_{k=1}^{+\infty} \int\frac{\1_{\{ N_\L(\om)=k \}}}{Z_{\L
,k}(\om_{\L^c})} \sum_{i=1}^k \int_{\L^k} (}{}\times \mathrm{e}^{-H_\L
(\om_{\L^c}\cup\listof x {1\ldots k} )}\dee\listelem x 1 \cdots\dee\listelem x k P(\mathrm{d}\om).
\nonumber
\end{eqnarray}
For every $\om$, by definition $(\nabla g+g\nabla H)(\xx,\om)$ is defined
for ${\lambda}^d$-almost every $\xx$, so the integral in
(\ref{devCamp})
is well-defined (either finite or infinite).

In the hereditary setting this equation characterizes
the canonical Gibbs measures if the energy functions $(H_\L)$
are sufficiently regular (Dereudre~\cite{DerCRASS} Theorem 1). However it is
clear that the characterization of canonical Gibbs measures via
(\ref{equationv}) in the non-hereditary setting is false in
general. Similar problems have been noted in Dereudre and Lavancier~\cite{DL};
see Remark~\ref{symetric} after Proposition~\ref{propvs} below.
\begin{pf*}{Proof of Proposition~\ref{propv}}
Let $g$ be a function in $\C_{\LL,\tau}$ satisfying (\ref{integrability})
with compact support contained in $\L=[-a,a]^d$.
The following calculations are available thanks to the
integrability assumption (\ref{integrability}).
As in (\ref{devCamp}), we have
%
\begin{equation}\label{pr1}
C_P^!(\nabla g - g\nabla H)
= \sum_{k=1}^{+\infty} \int\frac{\1_{\{ N_\L(\om)=k \}}}{Z_{\L
,k}(\om_{\L^c})}
\sum_{i=1}^k I^\L_{i,k}(\om_{\L^c}) P(\mathrm{d}\om),
\end{equation}
where
%
\begin{eqnarray}\label{int}\hspace*{-5pt}
I^\L_{i,k}(\om_{\L^c})
&=& \int_{\L^k} \bigl(
\nabla g(\listelem x i,\om_{\L^c}\cup\listof x {1\ldots k\bs i})\nonumber
\\[-12pt]
\\[-4pt]
&&\hphantom{\int_{\L^k} \bigl(}{} -g(\listelem x i,\om_{\L^c}\cup\listof x {1\ldots k\bs i})
\nabla H(\listelem x i,\om_{\L^c}\cup\listof x {1\ldots k\bs i})
\bigr)
\mathrm{e}^{-H_\L(\om_{\L^c}\cup\listof x {1\ldots k} )}
\dee\listelem x 1 \cdots\dee\listelem x k.\nonumber
\end{eqnarray}
Next, we show that $I^\L_{i,k}(\om_{\L^c})=0$ for every $\om_{\L^c}$,
every $k\ge1$ and every $1\le i \le k$. Consider the $j$th coordinate
of $I^\L_{i,k}(\om_{\L^c})$ which is denoted
by $\coord{ [I^\L_{i,k}(\om_{\L^c}) ]} j$.
By Fubini's theorem, and since $g$ is regularizing, we obtain
%
\begin{eqnarray}\label{intb}
&&\coord{ [ I^\L_{i,k}(\om_{\L^c}) ]} j\nonumber\\
&&\quad :=
\int_\L\ldots\int_\L \biggl(
\frac{\partial g}{\partial\coord x j}
(\listelem x i,\om_{\L^c}\cup\listof x {1\ldots k\bs i})- g(\listelem x i,\om_{\L^c}\cup\listof x {1\ldots k\bs i})
\frac{\partial H}{\partial\coord x j}
(\listelem x i,\om_{\L^c}\cup\listof x{1\ldots k\bs i})
\biggr)\nonumber\\
&&\quad\quad\hphantom{\int_\L\ldots\int_\L}{}\times \mathrm{e}^{-H_\L(\om_{\L^c}\cup\listelem x 1 \cup\cdots\cup\listelem x k )}
\dee\listelem x 1\cdots\dee\listelem x k
\nonumber\\
&&\quad = \int_\L\ldots\int_\L \biggl(
\int_{[-a,a]^{d-1}} \biggl(
\int_{-a}^a \frac{\partial}{\partial\coordlistelem x i j }
\bigl(
g(\listelem x i,\om_{\L^c}\cup\listof x {1\ldots k\bs i})
\mathrm{e}^{-H_\L(\listelem x i,\om_{\L^c}\cup\listof x {1\ldots k\bs i})}
\bigr)
\dee\coordlistelem x i j
\biggr)\nonumber
\\[-8pt]
\\[-8pt]
&&\quad\quad\hphantom{\int_\L\ldots\int_\L \biggl(
\int_{[-a,a]^{d-1}} }
\dee\coordlistelem x i 1 \cdots\dee\coordlistelem x i {j-1}
\dee\coordlistelem x i {j+1} \cdots\dee\coordlistelem x i d\biggr )
\dee\listelem x 1 \cdots\dee\listelem x {i-1} \dee\listelem x {i+1}
\dee\listelem x k
\nonumber\\
&&\quad= \int_\L\ldots\int_\L \biggl(
\int_{[-a,a]^{d-1}} \bigl[
g(\listelem x i,\om_{\L^c}\cup\listof x {1\ldots k\bs i})
\mathrm{e}^{-H_\L(\listelem x i,\om_{\L^c}\cup\listof x {1\ldots k\bs i})}
\bigr]^a_{\coordlistelem x i j =-a}
\nonumber\\
& &\quad\quad\hphantom{\int_\L\ldots\int_\L \biggl(
\int_{[-a,a]^{d-1}}} \dee\coordlistelem x i 1 \cdots\dee\coordlistelem x i {j-1}
\dee\coordlistelem x i {j+1}
\cdots\dee\coordlistelem x i d \biggr)
\dee\listelem x 1\cdots\dee\listelem x {i-1}
\dee\listelem x {i+1} \dee\listelem x k
\nonumber\\
&&\quad= 0.\nonumber
\end{eqnarray}
The last equality is due to the fact that $g$ has compact support
inside $\L$. This proves the proposition.
\end{pf*}
%
\subsection{Variational equation in the stationary case}
Here we give a counterpart of Proposition~\ref{propv} in the setting of
stationary point processes. As usual in this context, we assume the energy
functions $(H_\L)$ are \emph{shift-invariant}, meaning that for any vector
$u$ in $\R^d$, any configuration $\om$ and any $\L$, we have
$H_\L(\om)=H_{\tau_u(\L)}(\tau_u(\om))$.
\begin{proposition}\label{propvs}
Let $P$ be a stationary canonical Gibbs measure for
shift-invariant energy functions $(H_\L)$.
Write $\Eredpalm0$ for the expectation with respect to
the reduced Palm distribution $\redpalm0$ defined in (\ref{Palm}).
Then for every function $g$ in the class $\C_{\LL,\tau}$ such that
%
\begin{equation}\label{integrabilitys}
\Eredpalm0 \bigl(|\nabla g(0,\om)|+|g(0,\om)|+|g(0,\om)\nabla H(0,\om)|
\bigr)<\infty,
\end{equation}
we have
%
\begin{equation}\label{equationvs}
\Eredpalm0 (\nabla g(0,\om) )
= \Eredpalm0 (g(0,\om)\nabla H(0,\om) ).
\end{equation}
\end{proposition}

The most important difference between (\ref{equationv})
and (\ref{equationvs}) is that, in the latter, the function $g$ is
not assumed to have compact support with respect to the first
variable. This assumption was crucial to the proof of Proposition
\ref{propv}, as it ensures that the integral of the derivative is zero in
(\ref{intb}). Here stationarity replaces this assumption. This kind
of result was not observed in Almeida and Gidas~\cite{AG}. A condition like (1.2) in
Almeida and Gidas~\cite{AG} (namely $\int_{\R^n} \nabla\cdot (\mathbf{W}^{(\alpha)}\pi
_\theta(\xx) )\dee\xx=0$ for $\alpha=1,\ldots,m$)
is not required in this setting of stationary Gibbs point processes.
\begin{pf*}{Proof of Proposition~\ref{propvs}}
Let $g$ be a function in $\C_{\LL,\tau}$ satisfying
(\ref{integrabilitys}). For any $n\ge1$ we denote by $\L_n$ the set
$[-n,n]^d$ and $\partial\L_n$ the border of $\L_n$ defined by
$\L_n\bs\L_{n-1}$. Let $\psi_n$ be a differentiable partition of unity,
that is, a function from $\Rd$ to $\R$ such that $\psi$ is zero
outside $\L_n$
and equal to $1$ inside $\L_{n-1}$. We assume that $|\psi_n(\xx)|$ and
$|\nabla\psi_n(\xx)|$ are uniformly bounded with respect to $\xx$ and
$n$ by a finite positive constant $C_\psi$. Define the function
$g_n(\xx,\om)=\psi_n(\xx)g(\xx,\om)$. We claim that $g_n$
satisfies the assumptions of Proposition~\ref{propv}. Indeed it is
clear than $g_n$ has compact support, so $g_n$ is in
$\C_{\LL,K}$. It remains to check the integrability assumption
(\ref{integrabilitys}). From (\ref{Palm}) and noting that the functions
$g$, $\nabla g$ and $\nabla H$ are shift invariant, we have
%
\begin{eqnarray}\label{verif}
C_P^!(|\nabla g_n|+|g_n\nabla H|) & \le& z(P)\int\int (|\nabla\psi
_n||g|+|\psi_n||\nabla g|+|\psi_n|| g\nabla H| )\nonumber\\
&&\hphantom{z(P)\int\int}{}\times(\xx,\tau_\xx\om)
\redpalm0(\mathrm{d}\om){\lambda}^d(\mathrm{d}\xx)\\
& \le& z(P) {\lambda}^d(\L_n)C_\psi\int (|g|+|\nabla
g|+|g\nabla H| )(0,\om) \redpalm0 (\mathrm{d}\om)<\infty.\nonumber
\end{eqnarray}
Applying Proposition~\ref{propv} to $g_n$ we obtain
\[
\EE_P \biggl[
\sum_{\xx\in\L_n} \bigl(
\nabla g_n(\xx,\om\bs\xx)
- g_n(\xx,\om\bs\xx)\nabla H(\xx,\om\bs\xx)
\bigr) \biggr]=0.
\]
By stationarity, and since $g_n$ and $\nabla g_n$ on $\L_{n-1}$
are equivalent to $g$ and $\nabla g$, respectively, it follows that
%
\begin{eqnarray}\label{eqa}
0 & = & \EE_P \biggl[
\sum_{\xx\in\partial\L_n} \bigl(
\nabla\psi_n(\xx) g_n(\xx,\om\bs\xx)
+ \nabla g_n(\xx,\om\bs\xx)\psi_n(\xx)\nonumber\\
&&\hphantom{\EE_P \biggl[
\sum_{\xx\in\partial\L_n} \bigl(}{}- \psi_n(\xx,\om\bs\xx)g(\xx,\om\bs\xx)\nabla H(\xx,\om\bs
\xx)
\bigr) \biggr]\\
&&{} + z(P){\lambda}^d(\L_{n-1}) \Eredpalm0 \bigl(\nabla
g(0,\om)-g(0,\om)\nabla H(0,\om) \bigr).\nonumber
\end{eqnarray}
By a similar calculation as for (\ref{verif}), we deduce that the
first term
in (\ref{eqa}) is bounded by $K{\lambda}^d(\partial\L
_n)$ for $0 < K < \infty$.
Dividing (\ref{eqa}) by $n^d$ and letting $n$ tend to infinity,
the first term vanishes and we obtain (\ref{equationvs}).
\end{pf*}
\begin{remarque}\label{symetric}
For many choices of the function $g$
in Proposition~\ref{propvs}, the equation \eqref{equationvs}
is the trivial identity $0=0$.
For, suppose that the energy functions $(H_\L)$ are symmetric,
in the sense that $H_\L(\xx,\om)=H_{\vs(\L)}(\vs(\xx),\vs(\om))$
where $\vs$ is the symmetric transformation in $\R^d$
defined by $\vs(\xx)=-\xx$. Suppose that the measure $P$
is also symmetric, $P=P\circ\vs^{-1}$. This
situation applies to all examples in Section~\ref{SectionE}. Therefore,
for any symmetric function $g$ in $\C_{\LL,\tau}$, it follows that
$\nabla g$ and $\nabla H$ are anti-symmetric
(i.e., $\nabla g(\xx,\om)=-\nabla g(\vs(\xx),\vs(\om))$) and we deduce
\[
\Eredpalm0 (\nabla g(0,\om) )
= -\Eredpalm0 (\nabla g(0,-\om) )
= - \Eredpalm0 (\nabla g(0,\om) )\\
= 0.
\]
A similar calculation also gives $\Eredpalm0 (g(0,\om)\nabla H(0,\om
) )=0$.

This remark shows that, in order to obtain useful
instances of the equation (\ref{equationvs}),
the choice of function $g$ is delicate.
In the next section, we will see that an interesting choice
will be $g=\div H$ provided this belongs to $\C_{\LL,\tau}$.
In this situation, $g$ is anti-symmetric.
\end{remarque}
%
%
\section{Variational estimator procedure for exponential models}
\label{S:expo}
Now we assume the energy functions
depend on parameters $\tt=(\t_1,\ldots,\t_p)$ in $\R^p$
in the linear form
\[
H_\L^{\tt} = \tt\cdot\HH_\L:=\t_1 H^1_\L+\cdots+\t_p H^p_\L
\]
where $(\HH_\L) = ((H^1_\L),\ldots,(H^p_\L))$ is a suite of
$p$ families of energy functions. The resulting point process model is an
exponential family in the sense of Barndorff-Nielsen~\cite{barn78}, K{\"{u}}chler and S{{\o }}rensen~\cite{kuchsore97}.

In the remainder of the paper, we assume that $(\HH_\L)$ is shift invariant
and that, for some $\tt\in\R^p$, there exists a stationary canonical
Gibbs measure $P$
for the energy functions $(H_\L^{\tt})$.
In the following, we use $\C_{\LL,K}$ to denote the intersection of all
classes $\C_{\LL,K}((H^i_\L))$ for $1\le i \le p$, and similarly
for $\C_{\LL,\tau}$.

In this setting, we present two estimation procedures for $\tt$.
The first, called \emph{shift-invariant estimation} is based on the
equilibrium equation (\ref{equationvs}) in Proposition~\ref{propvs}.
This estimator is very natural and exploits the stationarity of the process.
The second procedure, called \emph{grid estimation}, is based on the
equilibrium equation (\ref{equationv}) in Proposition~\ref{propv}.
In this case, we subdivide the observation window using a grid,
and apply the equilibrium equation (\ref{equationv}) in each cell
of the grid. This procedure is less natural than the first,
but seems to enjoy better asymptotic properties.
In the next section, we show that both procedures are strongly consistent,
but we are only able to prove asymptotic normality for the second procedure.
%
\subsection{Shift-invariant estimation procedure}
From Proposition~\ref{propvs}, for any function $g$ in
the class $\C_{\LL,\tau}$ satisfying the integrability
assumption~(\ref{integrabilitys}) for any energy
functions $(H^{k}_\L)$, $1\le k \le p$, we obtain
%
\begin{equation}\label{equation}
\theta_1 \Eredpalm0 [ g(0,\om)\nabla H^1(0,\om) ]
+\cdots+
\theta_p \Eredpalm0 [ g(0,\om)\nabla H^p(0,\om) ]
= \Eredpalm0 [\nabla g(0,\om) ].
\end{equation}

This vectorial equation in dimension $d$ with $p$ parameters
gives a linear system of equations for $\tt$.
However, in many situations, the symmetry properties of the functions $g$
and the energy functions $(H^k_\L)$ (see examples in Section~\ref{SectionE})
give $d$ identical equations in (\ref{equation}).
Consequently, we keep only one equation in summing the $d$ equations
in (\ref{equation}). This seems to be the best strategy for extracting
maximum possible information from the data.
Therefore, in the following, the divergence operator
$\div=\frac{\partial}{\partial x^{(1)}}+\cdots+\frac{\partial
}{\partial x^{(d)}}$
is used in place of the gradient operator. Now\vspace*{1pt} a system
of $p$ equations is obtained by choosing $p$
functions $g_1,\ldots, g_p$ in $\C_{\LL,\tau}$
satisfying the integrability assumption (\ref{integrabilitys}).

Denoting by $A$ the $p \times p$ matrix
%
\begin{equation}\label{matrix}
A_{i,j}= \Eredpalm0 [g_i(0,\om)\div H^j(0,\om) ],
\end{equation}
and by $b=(b_1,\ldots,b_p)$ the $p$-vector
%
\begin{equation}\label{vector}
b_i= \Eredpalm0 [\div g_i(0,\om) ],
\end{equation}
the system of linear equations to determine $\tt$ is then $A\tt=b$.

Classically, we approximate $A$ and $b$ by the empirical average.
So for every $n\ge1$, denote by $\L_n$ the set $[-n,n]^d$
and for $P$-almost every realization $\om$ define the matrix
$\hat A^{(n)}$ with entries
\[
\hat A^{(n)}_{i,j}=\sum_{\xx\in\om_{\L_n}} g_i(\xx,\om)
\div H^j(\xx,\om),\vadjust{\goodbreak}
\]
and the vector $\hat b^{(n)}$ with entries
\[
\hat b^{(n)}_i=\sum_{\xx\in\om_{\L_n}}
\div g_i(\xx,\om).
\]
These are unnormalised empirical sums;
$(1/N_{\L_n}(\om)) \hat A^{(n)}$ and
$(1/N_{\L_n}(\om)) \hat b^{(n)}$ are consistent estimates of
$A$ and $b$, respectively.
If the matrix $\hat A^{(n)}$ is invertible,
we define the estimator $\hat\tt{}^{(n)}$ by
%
\begin{equation}\label{estimator}
\hat\tt{}^{(n)}= (\hat A^{(n)})^{-1} \hat b^{(n)}.
\end{equation}
Under suitable assumptions, the invertibility of $\hat A^{(n)}$
for sufficiently large $n$ and the strong consistency of this estimator
are proved in Proposition~\ref{propconvergence} of Section~\ref{SectionA}.
%
\subsection{Grid estimation procedure}
In this section, the cells of the grid are the cubes
$\D_u=\tau_u([0,1]^d)$ for any $u$ in $\Z^d$. By a classical rescaling
procedure, it is always possible to consider a grid with cubes
of side length $a>0$. For any function $g$ in
$\C_{\LL,K}$ with compact support in $\D_0$ and satisfying
integrability assumption (\ref{integrability}) for any energy
functions $(H^{i}_\L)$, $1\le i \le p$, we obtain by Proposition
\ref{propv}
%
\begin{eqnarray}\label{equationb}
&& E_{P}\biggl [\sum_{\xx\in\om_{\D_0}} \nabla g(\xx,\om\bs\xx)\biggr ] \nonumber\\
&&\quad = \theta_1 E_{P} \biggl[ \sum_{\xx\in\om_{\D_0}} g(\xx,\om\bs\xx
)\nabla H^1(\xx,\om\bs\xx) \biggr]+\cdots\\
&&\qquad{}+\theta_p E_{P} \biggl[\sum_{\xx
\in\om_{\D_0}} g(\xx,\om\bs\xx)\nabla H^p(\xx,\om\bs\xx)
\biggr].\nonumber
\end{eqnarray}
As in the first estimation procedure, we sum these equations and we use
the divergence operator in place of the gradient. So by choosing $p$
such functions $g_1,\ldots, g_p$,
we denote by $A$ the $p \times p$ matrix
%
\begin{equation}\label{matrixb}
A_{i,j}=E_{P} \biggl[ \sum_{\xx\in\om_{\D_0}} g_i(\xx,\om\bs\xx
)\div H^j(\xx,\om\bs\xx)\biggr ]
\end{equation}
and by $b=(b_1,\ldots,b_p)$ the $p$-vector
%
\begin{equation}\label{vectorb}
b_i=E_{P} \biggl[ \sum_{\xx\in\om_{\D_0}}\div g_i(\xx,\om\bs\xx) \biggr].
\end{equation}
As in the first estimation procedure, the system of linear equations to
determine $\tt$ is then $A\tt=b$.

For any function $g$ in $\C_{\LL,K}$ with compact support in $\D_0$,
we denote by $\bar g$ its periodic version defined by $\bar
g(\xx,\om)=\sum_{u\in\Z^d} g(\tau_u(\xx),\tau_u(\om))$. Then the
unnormalised empirical approximations of $A$ and $b$
are defined for every $n\ge1$
and $P$-almost every realization $\om$ by\looseness=-1
\begin{eqnarray*}
\hat A^{(n)}_{i,j} &=& \sum_{\xx\in\om_{\L_n}}
\bar g_i(\xx,\om\bs\xx)
\div H^j(\xx,\om\bs\xx),
\\
\hat b^{(n)}_i &=& \sum_{\xx\in\om_{\L_n}}
\div \bar g_i(\xx,\om\bs\xx).
\end{eqnarray*}\looseness=0
As in the first procedure, if the matrix $\hat A^{(n)}$ is invertible,
$\hat\tt{}^{(n)}$ is defined by
%
\begin{equation}\label{estimatorb}
\hat\tt{}^{(n)}= \bigl(\hat A^{(n)}\bigr)^{-1} \hat b^{(n)}.
\end{equation}

Under suitable assumptions, the invertibility of $\hat A^{(n)}$ for
sufficiently large $n$, strong consistency and asymptotic normality
of this estimator are proved in Propositions~\ref{propconvergence}
and~\ref{propconvergenceb} of Section~\ref{SectionA}.
%
%
\section{Asymptotic properties of the estimators}
\label{SectionA}
%
\subsection{Strong consistency}
\begin{proposition}\label{propconvergence}
Let $P$ be a stationary canonical Gibbs measure for shift invariant
energy functions $H_\L^{\tt}$. Let $g_1,\ldots,g_p$ be
functions in the class $\C_{\LL,\tau}$ (resp., in $\C_{\LL,K}$
with compact support in $\D_0$) satisfying the integrability assumption
(\ref{integrabilitys}) (resp., (\ref{integrability})) for all
$(H^i_\L)$, $1\le i \le p$. We assume also that, for any Gibbs measure
$\tilde P$ with energy functions $H_\L^{\tt}$, the matrix $A$ defined
in (\ref{matrix}) (resp., in (\ref{matrixb})) via $\tilde P$ is
invertible. Then for $P$-almost every realization $\om$, for sufficiently
large $n$ the matrix $\hat A^{(n)}$ is invertible and the estimator
$\hat
\tt{}^{(n)}$ defined in (\ref{estimator}) (resp., in
(\ref{estimatorb})) converges to $\tt$. That is, $\hat\tt{}^{(n)}$ is
strongly consistent.
\end{proposition}
\begin{pf}
Let $P$ be a stationary canonical Gibbs measure and $P = \int
P_z \nu(\mathrm{d} z)$ its mixture representation mentioned after
Definition~\ref{GCG}. Each Gibbs measure $P_z$ is itself a mixture of ergodic
Gibbs measures and therefore $P$ is also a mixture of ergodic Gibbs
measures (Georgii~\cite{G,geor88}).

So, for any $P$-a.e realization $\om$, $\om$ is also a realization of an
ergodic Gibbs measure~$\tilde P$. By ergodic theorem, the normalized
matrix $\hat A^{(n)}$ and vector $\hat b^{(n)}$ converge $\tilde
P$-a.s. to the matrix $A$ and vector $b$ defined via $\tilde
P$. Thanks to the equilibrium equations (\ref{equationv}),
(\ref{equationvs}) and the assumption that $A$ is invertible for any
Gibbs measure, the rest of the proof is clear.\looseness=-1
\end{pf}

As in Almeida and Gidas~\cite{AG}, it is not easy to find general conditions on
$g_1,\ldots,g_p$ which ensure that the matrix $A$ is invertible.
Nevertheless, there is one interesting choice for which it
is easy to prove it. Let us develop this situation in the rest of the
section.

In Almeida and Gidas~\cite{AG}, the authors propose to define
$g_i=\div H^i$.
We follow this idea with a small modification.
Indeed, in general, there is no reason why
$\div H^i$ should be regularizing, satisfy the
integrability assumption (\ref{integrabilitys}) or have compact
support in $\D_0$. So we propose the following variant. Let us choose
a fixed nonnegative function $\Psi$ from $\R^d\times\O$ to $\R$ and
define
%
\begin{equation}\label{gi}
g_i=\Psi\div H^i,\qquad 1\le i\le p.
\end{equation}

In the setting of shift invariant estimation, we will assume
that $\Psi$ is shift invariant. In the setting of grid estimation,
we will assume that $\Psi$ has compact support contained in~$\D_0$.
\begin{proposition}\label{propident}
Suppose that the functions $(g_i)_{1\le i\le p}$ defined as in (\ref{gi})
satisfy the following identifiability assumption: for any $X\in\R^p$,
%
\begin{equation}\label{identificationb}
\mbox{the function } \sum_{\xx\in\D_0} \Psi(\xx,\om\bs\xx)
|\transpose X \cdot \div\HH(\xx,\om\bs\xx) |=0, \qquad P\mbox{-a.s.}\quad\mbox{iff}\quad X=0.
\end{equation}
Then the associated matrix $A$ defined in (\ref{matrix}) or (\ref
{matrixb}) is invertible as soon as all the terms are integrable.
\end{proposition}

Let us remark that in the case of shift-invariant estimation, the
identification assumption (\ref{identificationb}) can be reformulated
by: for any $X\in\R^p$,
%
\begin{equation}\label{identification}
\mbox{the function } \Psi(0,\om) \transpose X \cdot \div\HH(0,\om)=0,\qquad
\redpalm0\mbox{-a.s.}\quad \mbox{iff}\quad  X=0.
\end{equation}
\begin{pf*}{Proof of Proposition~\ref{propident}}
In the shift invariant setting, from (\ref{Palm}) we observe that the
expressions for $A$ in (\ref{matrix}) and (\ref{matrixb}) are equivalent
up to a multiplicative scalar $z(P)$. We shall prove that $A$ defined in
(\ref{matrixb}) is invertible, by showing it is positive-definite.
Let $X$ be a vector in $\R^p$. We have
\begin{eqnarray*}
\transpose X A X & = & \transpose X \cdot E_{P}\biggl [ \sum_{\xx\in\om_{\D
_0}} \Psi(\xx,\om\bs\xx) \div\HH(\xx,\om\bs\xx)
\transpose{ \bigl(\div\HH(\xx,\om\bs\xx) \bigr)}
\biggr]\cdot X\\
& = & E_{P} \biggl[\sum_{\xx\in\om_{\D_0}} \Psi(\xx,\om\bs\xx)
\bigl(\transpose X \cdot\div\HH(\xx,\om\bs\xx)\bigr )^2 \biggr].
\end{eqnarray*}
Since $\Psi$ is nonnegative, this quantity is nonnegative and thanks
to the identification assumption~(\ref{identificationb}) it is
positive as soon as $X\neq0$.
\end{pf*}
\begin{remarque}
When $g_i=\psi\div H^i$, some terms in the matrix $A$ and the vector
$b$ can be simplified if the point process $P$ is symmetric in each
direction. Indeed, following the arguments in Remark~\ref{symetric},
it is easy to show that $\Eredpalm0( \frac{\partial H^i}{\partial
\coord x k}\frac{\partial H^j}{\partial\coord x l})=0$ and
$\Eredpalm0 (
\frac{\partial^2 H^i}{\partial\coord x k \partial\coord x l}
)=0$
as soon as $k\neq l$. Therefore, in the setting of shift invariant
estimator, the matrix $A$ and vector $b$ have the following simpler
expression,
%
\begin{equation}\label{matrixA}
A_{i,j}= \Eredpalm0 \Biggl[\Psi(0,\om) \Biggl(\sum_{k=1}^d \frac{\partial
H^i}{\partial\coord x k}(0,\om)\frac{\partial H^j}{\partial\coord x
k}(0,\om) \Biggr) \Biggr],
\end{equation}
and
%
\begin{equation}\label{vectorB}
b_i= \Eredpalm0 [\Psi(0,\om)\Delta H^i(0,\om)+\div\Psi(0,\om
)H^i(0,\om) ],
\end{equation}
where $\Delta$ denotes the classical Laplacian operator $\sum_{k=1}^d
\frac{\partial^2} {\partial^2 \coord x k}$. Obviously an analogue
simplification occurs for the grid estimator and these modifications
should be incorporated in the computation of the empirical matrix $\hat
A_n$ and vector $\hat b_n$.
\end{remarque}

Note that the variational estimator of Almeida and Gidas~\cite{AG} is an example of a
time-invariance estimator Baddeley~\cite{Baddeley2000}, that is, it can be derived
from properties of the infinitesimal generator of a certain diffusion.
In our case, the variational estimator can again be viewed as a
time-invariance estimator, associated with the diffusion with drift~$\nabla H$.
%
\subsection{Asymptotic normality}
We have seen above that shift invariant estimation seems more
natural than the grid estimation in the context of stationary
processes. Nevertheless, in this section we prove asymptotic normality
for the grid estimator, while we did not succeed in showing it for the
shift invariant estimator. However, in the simulations presented in
Section~\ref{SectionS}, we do not notice difference between the
asymptotic properties of these both estimators.

A function $g$ on $\R^d\times\O$ is said to have a
\emph{finite range} $R$ with $0 < R < \infty$ if for all $\xx$ and
$\om$
\[
g(\xx,\om)=g\bigl(\xx,\om_{B(\xx,R)}\bigr).
\]

\begin{proposition}\label{propconvergenceb}
Let $P$ be a stationary ergodic Gibbs measure for shift invariant
energy functions $H_\L^{\tt}$. Let $g_1,\ldots,g_p$ be
functions in the class in $\C_{\LL,K}$ satisfying the assumptions of
Proposition~\ref{propconvergence}.
Moreover, we assume that the functions $(g_i)$ and $(\nabla H^i)$ have
finite range $R > 0$, and that for every $1\le i,j \le p$
%
\begin{equation}\label{moment}
\EE_P \biggl( \biggl|\sum_{\xx\in\om_{\D_0}} \bigl(\nabla g_i(\xx,\om\bs\xx
)+g_i\nabla H^j(\xx,\om\bs_xx) \bigr) \biggr|^3\biggr )<\infty.
\end{equation}
Then the estimator $\hat\tt{}^{(n)}$ is asymptotically normal,
%
\begin{equation}\label{convdist}
{\lambda}^d(\L_n)^{1/2}\bigl(\hat\tt{}^{(n)}-\tt
\bigr)\Rightarrow\mathcal{N}(0,A^{-1}\Sigma A),
\end{equation}
where $\Rightarrow$ denotes convergence in distribution as $n \to
\infty$,
and $\Sigma$ is the matrix defined in~(\ref{Sigma}).
\end{proposition}

Let us notice that in Proposition~\ref{propconvergenceb}, we assume
that $P$ is ergodic which ensures, in general, that $P$ is an extremal
Gibbs measure. If it is not the case, then by the classical mixture
argument, we prove that the left term in (\ref{convdist}) converges in
distribution to a mixture of normal laws.
\begin{pf}
From the definition (\ref{estimatorb}) we have
\[
\hat A^{(n)}\bigl(\hat\tt{}^{(n)}-\tt\bigr)= \hat b^{(n)}-\hat A^{(n)}\tt.
\]
Since the assumptions of Proposition~\ref{propconvergence} hold and
since $P$ is ergodic,
the matrix ${\lambda}^d(\L_n)^{-1}\hat A^{(n)}$ converges
almost surely to $A$,
which is invertible.
The proof will be complete if we can show that the vector
$Z_n={\lambda}^d(\L_n)^{-1/2} (b^{(n)}-\hat
A^{(n)}\tt) $
converges in distribution to $\mathcal{N}(0,\Sigma)$.
Denote by $Z_n^i$ the $i$th coordinate of $Z_n$, we have
\[
Z_n^i= {\lambda}^d(\L_n)^{-1/2}\sum_{u\in\Z
^d\cap[-n,n-1]^d} Y_u^i,
\]
where
\[
Y_u^i=\sum_{\xx\in\om_{\D_u}} \Biggl(\div\bar g_i(\xx,\om\bs\xx)
-\sum_{j=1}^p \t_j \bar g_i(\xx,\om\bs\xx)\div H^j(\xx,\om\bs
\xx)\Biggr)
\]
where again $\D_u = \tau_u([0,1]^d)$.

We apply Theorem 2.1 in Jensen and K{\"u}nsch~\cite{JK} to obtain a Central Limit theorem
for $Z_n$. Let us check the three fundamental assumptions (\ref
{hyp1}), (\ref{hyp2}) and (\ref{hyp3}) below.
For every $u$ in $\Z^d$, consider the neighborhood of $u$ defined by
\[
\V_u= \{ v\in\Z^d, \mbox{ there exist } \xx\in\D_u, y\in\D_v
\mbox{ such that } |\xx-y|\le R \}.
\]
From the finite range of $(g_i)$ and $(\nabla H_i)$, it is easy to
remark that
%
\begin{equation}\label{hyp1}
Y_u^i(\om)= Y_u^i\biggl(\bigcup_{v\in\V_u} \om_{\D_v}\biggr).
\end{equation}
From stationarity of $P$, shift invariance of $\HH$, definition of
$\bar g_i$ and integrability assumption~(\ref{moment}), we obtain that
for every $u\in\Zd$ and every $1\le i \le p$
%
\begin{equation}\label{hyp2}
\EE_P(|Y_u^i|^3)= \EE_P(|Y_0^i|^3) < \infty.
\end{equation}
In Proposition~\ref{propv}, the equation (\ref{equationv}) remains valid
if $P$ is replaced by $P(\cdot|\om_{\L^c})$ where $\L$ is the
compact support
of $g$. For, in the proof we show that $I^\L_{i,k}(\om_{\L^c})=0$
for every $\om_{\L^c}$, every $k\ge1$ and every $1\le i \le k$.
Since the function $g_i$ has a compact support in $\D_0$,
we deduce that for every $1\le i \le p$ and $u\in\Zd$
%
\begin{equation}\label{hyp3}
\EE_P(Y_u^i|\om_{\D_u^c})=0.
\end{equation}
Applying Theorem 2.1 in Jensen and K{\"u}nsch~\cite{JK}, we conclude that $Z_n$ converges in
distribution to $\mathcal{N}(0,\Sigma)$ where $\Sigma$ is the
$p\times p$
matrix defined by
%
\begin{equation}\label{Sigma}
\Sigma_{i,j}= \EE_P \biggl(\sum_{v\in\V_0} Y_0^i Y_v^i \biggr).
\end{equation}
The proposition is proved.\vadjust{\goodbreak}
\end{pf}

To prove asymptotic normality for the shift invariant estimator, we
can apply the same argument as in the proof in Proposition
\ref{propconvergenceb}. Nevertheless it is not possible to show that
$Z_n$ satisfies a Central Limit theorem via Theorem 2.1 in Jensen and K{\"u}nsch~\cite{JK}
because the fundamental property (\ref{hyp3}) fails in this
situation. An alternative solution would be to substitute (\ref{hyp3})
by some mixing properties of the Gibbs measure $P$. We do not
investigate this solution here since we think that our estimators are
interesting in the setting of rigid point processes with strong
interaction. It is well known that mixing properties are not
established in this setting. Asymptotic normality of the shift invariant
estimator is not established.
%
%
\section{Examples}
\label{SectionE}
In this section, we present examples which are amenable to the
estimation procedures described in Section~\ref{S:expo}.
The first example is a model with unbounded and finite range pairwise potential
without hardcore part. The second involves a model of hardcore
spheres with interaction.

For the shift-invariant estimators, we will always choose the functions
$g_i$ with the form
%
\begin{equation}\label{gib}
g_i(\xx,\om)=\Psi(\xx,\om) \div H^i(\xx,\om),\qquad 1\le i\le p,
\end{equation}
where $\Psi$ is a shift invariant function from $\R^d\times\O$ to
$\R$ which we will determine for each model.

For the grid estimator, we will choose the functions $g_i$ with the form
%
\begin{equation}\label{gibb}
g_i(\xx,\om)=\psi(\xx) \Psi(\xx,\om) \div H^i(\xx,\om),\qquad 1\le
i\le p,
\end{equation}
where $\psi$ is a continuous and piecewise differentiable function
from $\R^d$ to $\R$ with a compact support exactly equals to
$[0,1]^d$ and such that $\nabla\psi$ is bounded. This function is
fixed for all the models.
An example of such a function $\psi$ is
%
\begin{equation}
\label{periodic}
\psi(\xx)=\1_{[0,1]^d}(\xx) \prod_{i=1}^d x^{(i)}\bigl(1 - x^{(i)}\bigr).
\end{equation}
%
\subsection{Pairwise interaction model}
In this section we study a general, unbounded, pairwise
potential with finite range. The infinite range case could be also
investigated, but it is more complicated to present: tempered
configurations have to be introduced, and the integrability
assumptions are much more difficult to obtain. Moreover in
statistical applications, the infinite range case has limited interest
because the observation window is typically bounded.

We assume also that the interaction does not include a hard core;
this setting is addressed in the next section.
Let $(\vp_i)_{1\le i \le p}$ be $p$ twice differentiable functions
from $]0,+\infty[$ to $\R\cup\{+\infty\}$ with compact support.
We denote by $R_0$ the common range of all potentials
$(\vp_i)$ (i.e., a real $R_0>0$ such that
every function $\vp_i$ is null on $[R_0,+\infty[$). The energy
functions $(H_\L^i)$ are defined by
%
\begin{equation}\label{energypairwise}
H^i_\L(\om)=\mathop{\mathop{\sum}_{
\{\xx,\yy\}\in\om}}_{\{\xx,\yy\} \cap\L\neq\varnothing}
 \vp_i \bigl((\xx-\yy)^2 \bigr) .
\end{equation}
The infinite sum in (\ref{energypairwise}) is well
defined for every $\om$ in $\O$ since the functions $(\vp_i)$ have compact
supports. The classical way to define a pairwise interaction is in
general via the quantities $\vp(|\xx-\yy|)$ but in our setting it is
simpler to use the equivalent form $\vp((\xx-\yy)^2 )$ since we
will compute derivatives with respect to the coordinates.

The global energy functions of the system are then defined by the
linear combination $H^{\tt}_\L=\t_1 H^1_\L+\cdots+\t_p H^p_\L$ for
$\tt\in\R^p$. Define the potential $\vp^{\tt}= \t_1\vp_1+\cdots
+\t_p \vp_p$. We assume that the potential $r\to\vp^{\tt}(r^2)$ is
superstable and lower regular (conditions (SS) and (LR) in
Ruelle~\cite{Ruelle70}, pages 131--132). These conditions assume that
double sums of pair potential terms $\sum_{i \neq j} \vp^{\tt}(\xx
_i, \xx_j)$
can be bounded from below by suitable expressions: the details are
not required for this paper. These conditions
ensure the existence of a stationary Gibbs measure $P$
(Ruelle~\cite{Ruelle70}, Theorem 5.8) for any intensity $z>0$ and provide also
the following property
(Ruelle~\cite{Ruelle70}, Corollary 2.9): for every $R>0$,
there exists $\alpha>0$ such that
%
\begin{equation}\label{estimate}
\EE_P (\mathrm{e}^{\alpha N^2_{B(0,R)}} )<0.
\end{equation}
We emphasise that the component functions
$\vp_i(r^2)$ are not assumed to be superstable and lower regular
for each $i$.

For this example, the function $\Psi$ in (\ref{gib}) is the constant
function equal to $1$ and so
$g_i=\div H^i$ for $1\le i\le p$.
Observe that the energy functions
$(H^i)_\L$ are hereditary, so the local energy $h^i(\xx,\om)$ exists
and is defined by
%
\begin{equation}\label{energielocale}
h^i(\xx,\om)= \sum_{\yy\in\om} \vp_i \bigl((\xx-\yy)^2 \bigr).
\end{equation}
For every $\om\in\O$, the local energy \eqref{energielocale}
is twice differentiable at every $\xx\notin\om$ and so
%
\begin{equation}\label{nablaH}
g_i(\xx,\om)=\div H^i(\xx,\om)= \div h^i (\xx,\om)= 2\sum_{\yy
\in\om} \Biggl(\vp'_i \bigl((\xx-\yy)^2 \bigr)\sum_{k=1}^d \bigl(\coord x k-\coord y
k\bigr) \Biggr),
\end{equation}
and
%
\begin{eqnarray}\label{deltaH}
\div g_i(\xx,\om)&=&\div\circ\div H^i (\xx,\om)\nonumber
\\[-8pt]
\\[-8pt]
&=& 2\sum_{\yy\in
\om} \Biggl( d\vp'_i \bigl((\xx-\yy)^2 \bigr)+2 \vp''_i \bigl((\xx-\yy)^2 \bigr) \Biggl(\sum
_{k=1}^d \bigl(\coord x k-\coord y k\bigr) \Biggr)^2 \Biggr).\nonumber
\end{eqnarray}

Let us give a collection of assumptions that will ensure the
shift invariant estimator $\hat\tt{}^{(n)}$ is strongly consistent.
\begin{definition}
The function $\vp^{\tt}$ satisfies property \textup{(LB)} if
it has a lower bound $C_\vp\in\R$,
\[
\inf_{r>0} \vp^{\tt}(r) \ge C_\vp.
\]
It satisfies the explosion control property \textup{(EX)} if
there exists $C_b>0$ such that for every $1\le i \le p$
\[
\sup_{r>0} |\vp'_i(r)|^2 \mathrm{e}^{-\vp^{\tt}(r)} \le C_b \quad\mbox{and}\quad
\sup_{r>0} |\vp''_i(r)| \mathrm{e}^{-\vp^{\tt}(r)} \le C_b.
\]
It satisfies the linear independence property \textup{(I)}
if the functions $(\vp'_i)_{1\le i \le p}$ are linearly independent
in the vectorial space of continuous functions from $]0,+\infty]\to\R$.
\end{definition}
\begin{proposition}\label{pairwiseconsistent}
Let $(\vp_i)$ be a family of potentials and $\tt\in\R^p$ such that
$\vp^{\tt}$ is superstable and lower regular, and such that
properties \textup{(LB)}, \textup{(EX)} and \textup{(I)} hold. Then for
any Gibbs measure associated to $\vp^{\tt}$,
the shift invariant estimator or the grid estimator $\hat\tt{}^{(n)}$
is strongly consistent.
\end{proposition}
\begin{pf}
Thanks to Propositions~\ref{propconvergence} and~\ref{propident},
it remains to check that $(g_i)$ are in $\C_{\LL,\tau}$ and that
assumptions (\ref{integrabilitys}) and (\ref{identification}) hold
for any Gibbs measure $P$.
The expression (\ref{nablaH}) shows that $g_i$ is in $\C_{\LL,\tau}$
for every $1\le i \le p$.

Concerning (\ref{identification}),
from (\ref{Palm}) and (\ref{NZess}) we deduce that
%
\begin{equation}\label{NZ}
\redpalm0(\mathrm{d}\om)=\frac{1}{z(P)}\mathrm{e}^{-h^{\tt}(0,\om)}P(\mathrm{d}\om).
\end{equation}

By the DLR equations (\ref{DLR}), it follows that $\redpalm0$ is locally
absolutely continuous with respect to the Poisson point Process
$\pi^z$ with positive density. We deduce that the probability under
$\redpalm0$
that the configuration $\om_{B(0,R_0)}$ is reduced to a point
$\{\xx\}$ is positive and that the law of this single point $\xx$ is
absolutely continuous with respect to the Lebesgue measure on
$B(0,R_0)$. In this situation, for any $X$ in $\Rd$ the left term in
(\ref{identification}) is nothing more than
%
\begin{equation}\label{identifi}
\Psi(0,\om) \transpose X \cdot\div\HH(0,\om)=-2 \Biggl(\sum_{k=1}^d \coord
x k \Biggr) ( \vp'_1 (x^2),\ldots, \vp'_p (x^2) )\cdot X=0.
\end{equation}
From assumption (I), this implies that
$X=0$ and (\ref{identification}) is proved.

It remains to check (\ref{integrabilitys}).
From expressions (\ref{nablaH}), (\ref{deltaH}) and formula (\ref
{NZ}), we have
\begin{eqnarray*}
& & \Eredpalm0 \bigl(|\nabla g_i(0,\om)|+|g_i(0,\om)|+|g_i(0,\om)\nabla
H^j(0,\om)| \bigr)\\
& &\quad\le K \Eredpalm0 \biggl[ \sum_{y\in\om_{B(0,R_0)}} \bigl(|\vp'_i(y^2)|+
|\vp''_i(y^2)| \bigr)+ \biggl( \sum_{y\in\om_{B(0,R_0)}} |\vp'_i(y^2)| \biggr) \biggl(
\sum_{y\in\om_{B(0,R_0)}} |\vp'_j(y^2)|\biggr ) \biggr]\\
& &\quad\le \frac{K}{z(P)} E_{P} \biggl[ \sum_{y\in\om_{B(0,R_0)}} \bigl(|\vp
'_i(y^2)|+|\vp''_i(y^2)|\bigr )\mathrm{e}^{-h^{\tt}(0,\om)}\\
&&\hphantom{\quad\le \frac{K}{z(P)} E_{P} \biggl[}{} +N_{B(0,R_0)} \mathrm{e}^{-h^{\tt}(0,\om)} \biggl( \biggl( \sum_{y\in\om_{B(0,R_0)}}
\vp'_i(y^2)^2 \biggr) + \biggl( \sum_{y\in\om_{B(0,R_0)}} \vp'_j(y^2)^2 \biggr) \biggr)\biggr ],
\end{eqnarray*}
where $K>0$ is a constant. Using (LB) and (EX), we obtain
%
\begin{eqnarray}\label{borne}
& & \Eredpalm0 \bigl(|\nabla g_i(0,\om)|+|g_i(0,\om)|+|g_i(0,\om)\nabla
H^j(0,\om)| \bigr)\nonumber\\
& &\quad\le \frac{K}{z(P)} E_{P} \biggl[ \sum_{y\in\om_{B(0,R_0)}} \bigl(1+\vp
'_i(y^2)^2+|\vp''_i(y^2)| \bigr)\mathrm{e}^{-\vp^{\tt}(y^2)}\mathrm{e}^{-h^{\tt}(0,\om\bs
y)}\nonumber
\\[-8pt]
\\[-8pt]
&&\hphantom{\quad\le \frac{K}{z(P)} E_{P} \biggl[}{}+ N_{B(0,R_0)}
\sum_{y\in\om
_{B(0,R_0)}} \bigl( \vp'_i(y^2)^2+\vp'_j(y^2)^2 \bigr) \mathrm{e}^{-\vp^{\tt}(y^2)}
\mathrm{e}^{-h^{\tt}(0,\om\bs y)} \biggr]\nonumber\\
&&\quad \le K E_{P} \bigl( (1+4 C_b)N^2_{B(0,R_0)} \mathrm{e}^{C_\vp N_{B(0,R_0)}}\bigr ).\nonumber
\end{eqnarray}
From (\ref{estimate}) this quantity is finite.
\end{pf}

For the grid estimator, we need a stronger assumption to obtain
asymptotic normality. Assumption (EX) is replaced by (EXb):
\[
\sup_{r>0} |\vp'_i(r)|^6 \mathrm{e}^{-\vp^{\tt}(r)} \le C_b \quad\mbox{and}\quad
\sup_{r>0} |\vp''_i(r)|^3 \mathrm{e}^{-\vp^{\tt}(r)} \le C_b.
\]

\begin{proposition}\label{pairwisenormality}
Let $(\vp_i)$ be a family of potentials and $\tt\in\R^p$ such that
$\vp^{\tt}$ is superstable and lower regular, and such that
properties \textup{(LB)}, \textup{(EXb)} and \textup{(I)} hold. Then
for any Gibbs measure associated to $\vp^{\tt}$,
the grid estimator $\hat\tt{}^{(n)}$ is strongly consistent and
asymptotically normal.
\end{proposition}
\begin{pf}
Thanks to Propositions~\ref{propconvergence},~\ref{propident}
and~\ref{propconvergenceb}, it remains to check that $(g_i)$ are in
$\C_{\LL,K}$ and that they satisfy (\ref{integrability}),
(\ref{moment}) and (\ref{identificationb}). Expression (\ref{nablaH})
and the assumptions on $\psi$ show that $g_i$ is in $\C_{\LL,K}$
for every $1\le i \le p$. Since the function $\psi$ in (\ref{gib}) is
bounded and has bounded gradient, by a calculation similar to that in
(\ref{borne}) and from \textup{(EXb)} we prove that
(\ref{integrability}) and (\ref{moment}) hold. The new assumption
\textup{(EXb)} plays a crucial role in the proof of
(\ref{moment}).
\end{pf}

An example of a potential $\vp^{\tt}$ satisfying
all assumptions (SS), (LR), (LB), (EXb)
and (I) is the
Lennard-Jones potential
%
\begin{equation}\label{LJ}
\vp^{\tt}(r) =
\1_{[0,R_0]}(r) \biggl(\theta_1 \frac{1}{r^6}+\theta_2 \frac{1}{r^3} \biggr),
\end{equation}
with $\theta_1>0$ and $\theta_2\in\R$. The Lennard-Jones potential
was first introduced in Lennard-Jones~\cite{lenn24}
with the parametrization $\theta_1=4\varepsilon\sigma^{12}$
and $\theta_2=-4\varepsilon\sigma^{6}$.
In Section~\ref{SectionS}, we show some simulations of this model
and investigate the estimation of $\varepsilon$ and $\sigma$.
The result shows that our estimator is very efficient and in particular,
when the model is very rigid, it seems better than the classical Likelihood
procedures.
%
\subsection{Interacting hard sphere model}
In this section, we consider the classical model of hard spheres with
pairwise interaction. We assume in this section that the pair
potential has finite interaction range,
is smooth, and does not explode (i.e., is bounded
near $0$). Exploding potentials were studied in the previous section;
a~mixture of both models could be investigated without additional difficulties.

Let us consider $R_0>r_0>0$ and $p$ functions $(\vp_i)_{1 \le i \le
p}$ from $\R^+$ to $\R\cup\{+\infty\}$ assumed to be twice
differentiable on
$[r_0,+\infty[$ with continuous second derivative. We assume that
$\vp_i$ is equal to infinity on $[0,r_0[$ and is null on
$[R_0,+\infty[$. The energy functions $(H_\L^i)$ are defined as in
(\ref{energypairwise}) and the definition of $H^{\tt}$ and
$\vp^{\tt}$ follows. In this setting, $\vp^{\tt}$ is necessarily
superstable and
lower regular, so a stationary Gibbs measure $P$ exists.

In this example, the energy functions are hereditary and so the
local energy $h^i(\xx,\om)$ exists and is defined in
(\ref{energielocale}). This local energy is differentiable for
${\lambda}^d$-almost every $\xx$ and so expression (\ref
{nablaH}) holds with
the convention $\vp'_i(r)=0$ if $r<r_0$.

It is clear that the function $\div H^i(\xx,\om)$ is not regularizing
because the function defined in
(\ref{functl}) is not continuous. Therefore, the choice
$g_i=\div H^i$ is not available in this
setting and an expression of type $g_i=\Psi\div
H^i$ is necessary. We define the function $\Psi$ by
%
\begin{equation}\label{psi}
\Psi(\xx,\om)=\1_{\Oi}(\om) \prod_{\yy\in\om} \chi_{r_0,r_1}
\bigl((\xx-\yy)^2 \bigr),
\end{equation}
where $\chi_{r_0,r_1}$ is the real function defined for any $r_1>r_0$ by
%
\begin{equation}
\label{chi}
\chi_{r_0,r_1}(r)= \cases{
0 &\quad if $r\le r_0$, \cr
\displaystyle\frac{r-r_0}{r_1-r_0} & \quad if $r_0 \le r \le r_1$, \cr
1 &\quad if $r\ge r_1$.}
\end{equation}
The product in (\ref{psi}) in fact involves only a finite number of terms
and so is well defined.
It is clear that the function $\Psi$ is ${\lambda
}^d$-a.e. differentiable with
\begin{eqnarray}\label{psideriv}
\div\Psi(\xx,\om)&=&\frac{2}{r_1-r_0}\1_{\Oi}(\om)\nonumber
\\[-8pt]
\\[-8pt]
&&{}\times\sum_{\yy\in
\om} \Biggl(\1_{[r_0,r_1]} \bigl((\xx-\yy)^2 \bigr) \Biggl(\sum_{k=1}^d \coord x k
-\coord y k \Biggr)\prod_{\zz\in\om\bs\yy} \chi_{r_0,r_1} \bigl((\xx-\zz
)^2 \bigr) \Biggr).\nonumber
\end{eqnarray}

\begin{proposition}\label{prophardcore}
If the potentials $(\vp'_i)$ are linearly independent
(assumption \textup{(I)}) then the shift invariant estimator and
the grid estimator $\hat\tt{}^{(n)}$ are strongly consistent.
Moreover, the grid estimator is asymptotically normal.
\end{proposition}
\begin{pf}
Thanks to Propositions~\ref{propconvergence},~\ref{propident} and
\ref{propconvergenceb}, it remains to check that $(g_i)$ are in $\C
_{\LL,\tau}$ and that they satisfy (\ref{integrability}), (\ref
{identificationb}) and (\ref{moment}).
First, from expression (\ref{nablaH}) and (\ref{psi}) let us note
that $g_i$ is regularising since the function defined in (\ref
{functl}) is always continuous and piecewise differentiable with
continuous derivatives. Since $g_i$ is clearly shift invariant it
follows that $g_i$ is in $\C_{\LL,\tau}$. For integrability
assumptions (\ref{integrability}), (\ref{moment}), it is sufficient
to remark that $g_i$, $\nabla g_i$ and $g_i\nabla H_j$ are uniformly
bounded. Finally for identification assumption (\ref
{identificationb}), the proof is exactly the same as in Proposition
\ref{pairwiseconsistent}.
\end{pf}

Let us remark that, in the case of hard packing (i.e., when $z$ is very large),
the Maximum Likelihood and Maximum Pseudolikelihood procedures are not
really available to estimate $\tt$. Our procedure should be efficient.

The variational estimator can also be applied to models with
rigid geometric constraints, such as the rigid Voronoi models of Dereudre and Lavancier~\cite{DL}.

\section{Simulations}
\label{SectionS}
In this section, we present the results of
simulation experiments assessing the performance of parameter estimation
in the Lennard-Jones model (\ref{LJ}). We consider four cases
where $\varepsilon$ takes the values $0.1$, $0.5$, $1$ and $2$, respectively,
which we call the cases of low, moderate, high and extreme rigidity.
We chose intensity $z = 100$
and characteristic range $\sigma=0.1$. For each model, we simulated
$1000$ realizations in the window $[0,2]^2$ and estimated the
parameters $\varepsilon$ and $\sigma$ (via $\theta_1=4\varepsilon
\sigma^{12}$ and $\theta_2=-4\varepsilon\sigma^{6}$) using the shift
invariant estimator \eqref{estimator}
and the grid estimator \eqref{estimatorb} where the window was
subdivided into a $10\times10$ grid of squares.
The terms in \eqref{estimator} were computed using
\eqref{nablaH} and \eqref{deltaH} with
$\vp_1(s) = s^{-6}$ and $\vp_2(s) = s^{-3}$. Similarly for \eqref{estimatorb}
using the periodic function $\psi$ in \eqref{periodic}.
For comparison, we also computed the maximum pseudolikelihood estimator (MPLE)
using the method of Baddeley and Turner~\cite{baddturn00}.
%

\begin{table}
\tablewidth=\textwidth
\tabcolsep=0pt
\caption{Sample statistics for 1000 estimates
of the parameters $\varepsilon, \sigma$ for the Lennard-Jones model
using the variational estimators (\textsc{grid} and \textsc{invariant})
and the approximate maximum pseudolikelihood estimator (\textsc{mpl})}\label{tab:EstiLJ}
\begin{tabular*}{\textwidth}{@{\extracolsep{\fill}}llld{4.3}llllll@{}}
\hline
&   &
\multicolumn{2}{l}{Low} &
\multicolumn{2}{l}{Moderate} &
\multicolumn{2}{l}{High} &
\multicolumn{2}{l}{Extreme} \\[-5pt]
&   &
\multicolumn{2}{l}{\hrulefill} &
\multicolumn{2}{l}{\hrulefill} &
\multicolumn{2}{l}{\hrulefill} &
\multicolumn{2}{l@{}}{\hrulefill} \\
&  & $\sigma$ & \multicolumn{1}{l}{$\varepsilon$} & $\sigma$ &
$\varepsilon$ &
$\sigma$ & $\varepsilon$ & $\sigma$ &
$\varepsilon$ \\
\hline
\textsc{true} & &
0.1 & 0.1 & 0.1 & 0.5
& 0.1 & 1 & 0.1 & 2 \\[3pt]
\textsc{grid}
& mean & 0.094 & 0.789 & 0.101 & 0.966 & 0.101 & 1.277 & 0.100 & 1.987
\\
& median & 0.091 & 0.570 & 0.098 & 0.882 & 0.099 & 1.244 & 0.100 &
1.946 \\
& sd & 0.016 & 2.253 & 0.011 & 0.715 & 0.007 & 0.691 & 0.003 & 0.822 \\
& \textsf{E}& 0.099 & 0.124 & 0.100 & 0.511 & 0.100 & 1.019 & 0.100 &
1.807 \\[3pt]
\textsc{invariant}
& mean & 0.093 & 0.636 & 0.100 & 0.970 & 0.100 & 1.333 & 0.099 & 2.210
\\
& median & 0.091 & 0.576 & 0.098 & 0.870 & 0.099 & 1.300 & 0.099 &
2.164 \\
& sd & 0.013 & 4.491 & 0.010 & 0.663 & 0.006 & 0.653 & 0.003 & 0.704 \\
& \textsf{E} & 0.097 & 0.149 & 0.099 & 0.558 & 0.099 & 1.106 & 0.099 &
2.069 \\[3pt]
\textsc{mpl}
& mean & 0.068 & -614.696 & 0.102 & 0.311 & 0.102 & 0.307 & 0.103 &
0.327 \\
& median & 0.094 & 0.050 & 0.101 & 0.314 & 0.102 & 0.299 & 0.103 &
0.273 \\
& sd & 0.040 & 763.422 & 0.004 & 0.129 & 0.002 & 0.096 & 0.002 & 0.168\\
\hline
\end{tabular*}  \vspace*{-3pt}
\end{table}

\begin{table}[b]\vspace*{-3pt}
\tablewidth=\textwidth
\tabcolsep=0pt
\caption{Summary of estimates of canonical
parameters $\theta_1, \theta_2$ in the same experiments
as previous table.
Values of $\theta_1$ multiplied by $10^{12}$;
values of $\theta_2$ multiplied by $10^6$}
\label{tab:EstiLJ:theta}
\begin{tabular*}{\textwidth}{@{\extracolsep{\fill}}llld{2.3}ld{2.3}ld{2.3}ld{2.3}@{}}
\hline
& &
\multicolumn{2}{l}{Low} &
\multicolumn{2}{l}{Moderate} &
\multicolumn{2}{l}{High} &
\multicolumn{2}{l}{Extreme} \\[-5pt]
& &
\multicolumn{2}{l}{\hrulefill} &
\multicolumn{2}{l}{\hrulefill} &
\multicolumn{2}{l}{\hrulefill} &
\multicolumn{2}{l@{}}{\hrulefill} \\
&&
$\theta_1$ & \multicolumn{1}{l}{$\theta_2$} & $\theta_1$ & \multicolumn{1}{l}{$\theta_2$} &
$\theta_1$ & \multicolumn{1}{l}{$\theta_2$} & $\theta
_1$ & \multicolumn{1}{l}{$\theta_2$}\\
 \hline
\textsc{true} & &
0.4 & -0.4 & 2 & -2 &
4 & -4 & 8 & -8 \\[3pt]
\textsc{grid}
& mean & 0.777 & -1.371 & 2.705 & -3.107 & 4.496 & -4.718 & 7.462 &
-7.636 \\
& median & 0.671 & -1.163 & 2.552 & -2.990 & 4.392 & -4.695 & 7.364 &
-7.538 \\
& sd & 0.511 & 1.296 & 1.210 & 1.945 & 1.370 & 2.036 & 1.553 & 2.387 \\
& \textsf{E}& 0.416 & -0.454 & 2.018 & -2.031 & 4.020 & -4.048 & 7.123
& -7.176 \\  [3pt]
\textsc{invariant}
& mean & 0.758 & -1.362 & 2.682 & -3.127 & 4.535 & -4.86 & 7.668 &
-8.199 \\
& median & 0.663 & -1.212 & 2.540 & -2.947 & 4.463 & -4.814 & 7.608 &
-8.127 \\
& sd & 0.470 & 1.178 & 1.147 & 1.805 & 1.290 & 1.911 & 1.379 & 2.033 \\
& \textsf{E}& 0.421 & -0.502 & 2.057 & -2.143 & 4.116 & -4.267 & 7.386
& -7.819 \\[3pt]
\textsc{mpl}
& mean & 0.274 & -0.096 & 1.435 & -1.317 & 1.526 & -1.360 & 1.715 &
-1.487 \\
& median & 0.291 & -0.240 & 1.410 & -1.332 & 1.509 & -1.355 & 1.602 &
-1.309 \\
& sd & 0.265 & 0.575 & 0.350 & 0.419 & 0.266 & 0.316 & 0.471 & 0.557 \\
\hline
\end{tabular*}
\end{table}

Table~\ref{tab:EstiLJ} summarises the 1000 estimated values for each
experiment. For each estimation method, we report the sample mean,
median and standard deviation of the individual estimates $\hat\sigma$
and $\hat\varepsilon$ for each of the three estimation methods. As a
cross-check on the validity of the method, the row marked \textsf{E} in
Table~\ref{tab:EstiLJ} denotes the `pooled' estimate of the parameters
obtained by computing the averages $A^\ast, b^\ast$ of the 1000
empirical matrices and vectors $\hat A, \hat b$, then solving
$A^\ast\tt= b^\ast$ and transforming to obtain $(\sigma, \varepsilon)$.
Table~\ref{tab:EstiLJ:theta} summarises the corresponding estimates of
the canonical parameters $\theta_1,\theta_2$ for the same experiment.

\begin{figure}[t!]

\includegraphics{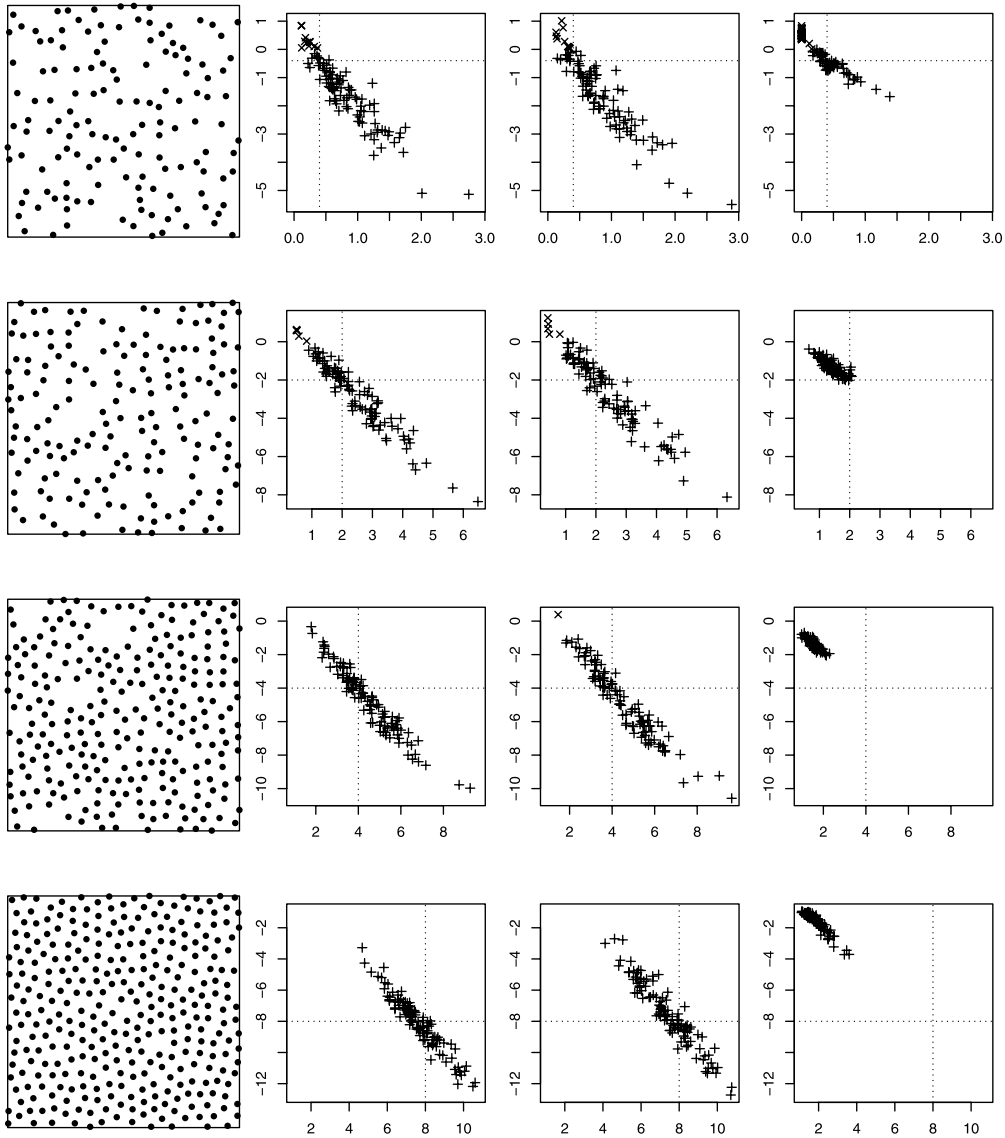}

\caption{Distribution of parameter estimates
for the Lennard-Jones model.
Each row corresponds to a different experiment, obtained by simulating
the Lennard-Jones model with (from top to bottom)
$\varepsilon=0.1, 0.5, 1, 2$, respectively.
The left column shows a typical realization.
The remaining columns show estimates of $(\theta_1,\theta_2)$
from 100 realizations, obtained using (from left to right)
the shift invariant estimator,
grid estimator and maximum pseudolikelihood estimator respectively.
Dotted lines show the true parameter values.}\label{fig:EstiLJ}
\end{figure}

Figure~\ref{fig:EstiLJ} displays typical realisations of the model
in each case, and a scatterplot of a subsample of 100 estimated values
giving an impression of the probability distribution of the estimates.

The most striking feature of these simulations is the bias in the
maximum pseudolikelihood estimator. The bias increases with
$\varepsilon$. We computed the MPLE using the Berman--Turner device
Baddeley and Turner~\cite{baddturn00}, Berman and Turner~\cite{bermturn92} with a $256 \times256$ grid of sample
points. To minimise numerical problems (overflow, instability, slow
convergence), we rescaled the interpoint distances to a unit equal to
the true value of $\sigma$. This is unrealistic with respect to
applications (since the true $\sigma$ would not be known) but gives
the most optimistic assessment of performance from the MPLE algorithm.
Although the MPLE is known to be biased in the presence of strong interaction,
we conjecture that the very large bias observed here may be due to
discretization error (Baddeley and Turner~\cite{baddturn00}, equation (17), page 290) and
numerical problems related to the pair potential.

The parameter space for the canonical parameters is
$\Theta= \{ (\theta_1,\theta_2)\dvt \theta_1 > 0, \theta_2 < 0\}$.
Our implementation of the three parameter estimation methods
is unconstrained, so that the algorithms may yield estimates
$(\hat\theta_1, \hat\theta_2)$ that are ``invalid'' in the sense
that they lie outside~$\Theta$. Invalid estimates are plotted as
crosses in Figure~\ref{fig:EstiLJ}. The observed proportion of
invalid estimates is shown in Table~\ref{tab:invalid}. In the case
of an invalid estimate, the model should be refitted
with $\theta$ constrained to lie on the boundary of $\Theta$
(i.e., constraining either $\theta_1 = 0$ or $\theta_2 = 0$).

The variational estimator seems to work
better when the model is rigid. In the low rigidity case, the
estimation of $\varepsilon$ is severely biased.
We interpret the bias as a consequence of
the nonlinear relationship between $\tt$ and $(\sigma, \varepsilon)$
combined with high variability in the low rigidity case.

The standard deviation of $\hat\varepsilon$ remains high in all cases.
This is unsurprising since the estimation of
$\varepsilon$ is a difficult problem for all
standard procedures even if the model is not rigid.\looseness=-1

The difference between the shift invariant
estimator and the grid estimator seems rather slight overall.
The better asymptotic
properties proved for the grid estimator seem to be irrelevant in these
simulations.
The estimates of $\sigma$ and $\varepsilon$ are strongly negatively correlated
in all cases. In the low rigidity case, the empirical distribution
of $\hat\varepsilon$ appears to deviate substantially from a Normal distribution,
suggesting that the behaviour predicted by the Central Limit theorem
has not yet set in.

However, the \emph{pooled} grid estimate of $\varepsilon$
appears to be biased in the extremely rigid case.
We attribute this to the absence of edge correction in our algorithm.

\begin{table}
\tablewidth=\textwidth
\tabcolsep=0pt
\caption{Proportion of estimates that were invalid,
in each of the four experiments}
\label{tab:invalid}
\begin{tabular*}{\textwidth}{@{\extracolsep{\fill}}lllll@{}}
\hline
& Low & Moderate & High & Extreme \\
\hline
\textsc{grid} & 0.116 & 0.031 & 0.006 & 0.000 \\
\textsc{Invariant} & 0.105 & 0.023 & 0.004 & 0.000 \\
\textsc{mpl} & 0.413 & 0.001 & 0.000 & 0.000 \\
\hline
\end{tabular*} \vspace*{-3pt}
\end{table}

In these experiments, the maximum pseudolikelihood estimator
had substantially greater computational cost than the variational estimators.
All algorithms were implemented in the \textsf{R} language using the spatial
statistics package \texttt{spatstat} Baddeley and Turner~\cite{baddturn05}.
Average computation times on a 2.5 GHz laptop were
about 7 seconds to generate one simulated realisation,
about 0.5 seconds each for the variational estimators,
and for the MPLE, about 2 sec, 4 sec, 7 sec and 60 sec
for the low, medium, high and extreme rigidity cases respectively.
A supercomputing cluster was used to conduct the simulations
in the extreme-rigidity case.

\section*{Acknowledgements}
We thank Mark Westcott and two anonymous referees for insightful comments.
Lawrence Murray and John Taylor provided invaluable advice on supercomputing.


\printhistory


\begin{thebibliography}{35}

\bibitem{AG}
\begin{barticle}[mr]
\bauthor{\bsnm{Almeida},~\bfnm{Murilo~P.}\binits{M.P.}} \AND
  \bauthor{\bsnm{Gidas},~\bfnm{Basilis}\binits{B.}}
(\byear{1993}).
\btitle{A variational method for estimating the parameters of {MRF} from
  complete or incomplete data}.
\bjournal{Ann. Appl. Probab.}
\bvolume{3}
\bpages{103--136}.%
\bid{issn={1050-5164}, mr={1202518}}%
\bptok{imsref}%
\end{barticle}%
\endbibitem\vadjust{\goodbreak}

\bibitem{baddturn00}
\begin{barticle}[mr]
\bauthor{\bsnm{Baddeley},~\bfnm{Adrian}\binits{A.}} \AND
  \bauthor{\bsnm{Turner},~\bfnm{Rolf}\binits{R.}}
(\byear{2000}).
\btitle{Practical maximum pseudolikelihood for spatial point patterns (with
  discussion)}.
\bjournal{Aust. N. Z. J. Stat.}
\bvolume{42}
\bpages{283--322}.
\bid{doi={10.1111/1467-842X.00128}, issn={1369-1473}, mr={1794056}}
\bptok{imsref}%
\end{barticle}
\endbibitem

\bibitem{baddturn05}
\begin{barticle}[auto:STB|2012/03/20|13:03:14]
\bauthor{\bsnm{Baddeley},~\bfnm{A.}\binits{A.}} \AND
  \bauthor{\bsnm{Turner},~\bfnm{R.}\binits{R.}}
(\byear{2005}).
\btitle{Spatstat: An {\textsf{R}} package for analyzing spatial point
  patterns}.
\bjournal{J.~Statist. Softw.}
\bvolume{12}
\bpages{1--42}.
\bptok{imsref}%
\end{barticle}
\endbibitem

\bibitem{Baddeley2000}
\begin{barticle}[mr]
\bauthor{\bsnm{Baddeley},~\bfnm{A.~J.}\binits{A.J.}}
(\byear{2000}).
\btitle{Time-invariance estimating equations}.
\bjournal{Bernoulli}
\bvolume{6}
\bpages{783--808}.
\bid{doi={10.2307/3318756}, issn={1350-7265}, mr={1791902}}
\bptok{imsref}%
\end{barticle}
\endbibitem

\bibitem{barn78}
\begin{bbook}[mr]
\bauthor{\bsnm{Barndorff-Nielsen},~\bfnm{Ole}\binits{O.}}
(\byear{1978}).
\btitle{Information and Exponential Families in Statistical Theory}.
\baddress{Chichester}: \bpublisher{Wiley}.
\bid{mr={0489333}}
\bptok{imsref}%
\end{bbook}
\endbibitem

\bibitem{bermturn92}
\begin{barticle}[auto:STB|2012/03/20|13:03:14]
\bauthor{\bsnm{Berman},~\bfnm{M.}\binits{M.}} \AND
  \bauthor{\bsnm{Turner},~\bfnm{T.~R.}\binits{T.R.}}
(\byear{1992}).
\btitle{Approximating point process likelihoods with GLIM}.
\bjournal{Applied Statistics}
\bvolume{41}
\bpages{31--38}.
\bptok{imsref}%
\end{barticle}
\endbibitem

\bibitem{besa78}
\begin{barticle}[mr]
\bauthor{\bsnm{Besag},~\bfnm{Julian}\binits{J.}}
(\byear{1977}).
\btitle{Some methods of statistical analysis for spatial data}.
\bjournal{Bull. Inst. Internat. Statist.}
\bvolume{47}
\bpages{77--91, 138--147}.
\bid{mr={0617572}}
\bptnote{check related, check year}
\bptok{imsref}%
\end{barticle}
\endbibitem



\bibitem{billcoeudrou08}
\begin{barticle}[mr]
\bauthor{\bsnm{Billiot},~\bfnm{Jean-Michel}\binits{J.M.}},
  \bauthor{\bsnm{Coeurjolly},~\bfnm{Jean-Fran{\c{c}}ois}\binits{J.F.}} \AND
  \bauthor{\bsnm{Drouilhet},~\bfnm{R{\'e}my}\binits{R.}}
(\byear{2008}).
\btitle{Maximum pseudolikelihood estimator for exponential family models of
  marked {G}ibbs point processes}.
\bjournal{Electron. J.~Stat.}
\bvolume{2}
\bpages{234--264}.
\bid{doi={10.1214/07-EJS160}, issn={1935-7524}, mr={2399195}}
\bptok{imsref}%
\end{barticle}
\endbibitem

\bibitem{CDDL}
\begin{bmisc}[auto:STB|2012/03/20|13:03:14]
\bauthor{\bsnm{Coeurjolly},~\bfnm{J.~F.}\binits{J.F.}},
  \bauthor{\bsnm{Dereudre},~\bfnm{D.}\binits{D.}},
  \bauthor{\bsnm{Drouilhet},~\bfnm{R.}\binits{R.}} \AND
  \bauthor{\bsnm{Lavancier},~\bfnm{F.}\binits{F.}}
(\byear{2010}).
\bhowpublished{Takacs--{F}iksel method for stationary marked Gibbs point
  processes. Preprint}.
\bptok{imsref}%
\end{bmisc}
\endbibitem

\bibitem{coeudrou10}
\begin{barticle}[mr]
\bauthor{\bsnm{Coeurjolly},~\bfnm{Jean-Fran{\c{c}}ois}\binits{J.F.}} \AND
  \bauthor{\bsnm{Drouilhet},~\bfnm{R{\'e}my}\binits{R.}}
(\byear{2010}).
\btitle{Asymptotic properties of the maximum pseudo-likelihood estimator for
  stationary {G}ibbs point processes including the {L}ennard-{J}ones model}.
\bjournal{Electron. J. Stat.}
\bvolume{4}
\bpages{677--706}.
\bid{doi={10.1214/09-EJS494}, issn={1935-7524}, mr={2678967}}
\bptok{imsref}%
\end{barticle}
\endbibitem


\bibitem{DerCRASS}
\begin{barticle}[mr]
\bauthor{\bsnm{Dereudre},~\bfnm{David}\binits{D.}}
(\byear{2002}).
\btitle{Une caract\'erisation de champs de {G}ibbs canoniques sur {${\Bbb R}\sp
  d$} et {${\mathcal C}([0,1],{\Bbb R}\sp d)$}}.
\bjournal{C.~R. Math. Acad. Sci. Paris}
\bvolume{335}
\bpages{177--182}.
\bid{doi={10.1016/S1631-073X(02)02452-4}, issn={1631-073X}, mr={1920016}}
\bptok{imsref}%
\end{barticle}
\endbibitem

\bibitem{Dereudre}
\begin{barticle}[mr]
\bauthor{\bsnm{Dereudre},~\bfnm{David}\binits{D.}}
(\byear{2008}).
\btitle{Gibbs {D}elaunay tessellations with geometric hardcore conditions}.
\bjournal{J. Stat. Phys.}
\bvolume{131}
\bpages{127--151}.
\bid{doi={10.1007/s10955-007-9479-6}, issn={0022-4715}, mr={2394701}}
\bptok{imsref}%
\end{barticle}
\endbibitem

\bibitem{DDG}
\begin{barticle}[auto:STB|2012/03/20|13:03:14]
\bauthor{\bsnm{Dereudre},~\bfnm{D.}\binits{D.}},
  \bauthor{\bsnm{Drouilhet},~\bfnm{R.}\binits{R.}} \AND
  \bauthor{\bsnm{Georgii},~\bfnm{H.~O.}\binits{H.O.}}
(\byear{2012}).
\btitle{Existence of Gibbsian point processes with geometry-dependent
  interactions}.
\bjournal{Probab. Theory Related Fields}
\bvolume{153}
\bpages{643--670}.
\bid{doi={10.1007/s00440-011-0356-5}, mr={2948688}}
\bptok{imsref}%
\end{barticle}
\endbibitem

\bibitem{DL2}
\begin{barticle}[auto:STB|2012/03/20|13:03:14]
\bauthor{\bsnm{Dereudre},~\bfnm{D.}\binits{D.}} \AND
  \bauthor{\bsnm{Lavancier},~\bfnm{F.}\binits{F.}}
(\byear{2011}).
\btitle{Practical simulation and estimation for Gibbs Delaunay--{V}oronoi
  tessellations with geometric hardcore interaction}.
\bjournal{Comput. Statist. Data Anal.}
\bvolume{55}
\bpages{498--519}.
\bid{mr={2736572}}
\bptok{imsref}%
\end{barticle}
\endbibitem

\bibitem{DL}
\begin{barticle}[mr]
\bauthor{\bsnm{Dereudre},~\bfnm{David}\binits{D.}} \AND
  \bauthor{\bsnm{Lavancier},~\bfnm{Fr{\'e}d{\'e}ric}\binits{F.}}
(\byear{2009}).
\btitle{Campbell equilibrium equation and pseudo-likelihood estimation for
  non-hereditary {G}ibbs point processes}.
\bjournal{Bernoulli}
\bvolume{15}
\bpages{1368--1396}.
\bid{doi={10.3150/09-BEJ198}, issn={1350-7265}, mr={2597597}}
\bptok{imsref}%
\end{barticle}
\endbibitem

\bibitem{fiks84}
\begin{barticle}[mr]
\bauthor{\bsnm{Fiksel},~\bfnm{Thomas}\binits{T.}}
(\byear{1984}).
\btitle{Estimation of parametrized pair potentials of marked and nonmarked
  {G}ibbsian point processes}.
\bjournal{Elektron. Informationsverarb. Kybernet.}
\bvolume{20}
\bpages{270--278}.
\bid{issn={0013-5712}, mr={0760146}}
\bptok{imsref}%
\end{barticle}
\endbibitem

\bibitem{G}
\begin{bbook}[mr]
\bauthor{\bsnm{Georgii},~\bfnm{Hans-Otto}\binits{H.O.}}
(\byear{1979}).
\btitle{Canonical {G}ibbs Measures: Some Extensions of de Finetti's Representation Theorem for Interacting
  Particle Systems}.
\bseries{Lecture Notes in Math.}
\bvolume{760}.
\baddress{Berlin}: \bpublisher{Springer}.
\bid{mr={0551621}}
\bptok{imsref}%
\end{bbook}
\endbibitem

\bibitem{geor88}
\begin{bbook}[mr]
\bauthor{\bsnm{Georgii},~\bfnm{Hans-Otto}\binits{H.O.}}
(\byear{1988}).
\btitle{Gibbs Measures and Phase Transitions}.
\bseries{de Gruyter Studies in Mathematics}
\bvolume{9}.
\baddress{Berlin}: \bpublisher{de Gruyter}.
\bid{mr={0956646}}
\bptok{imsref}%
\end{bbook}
\endbibitem

\bibitem{geye99}
\begin{bincollection}[mr]
\bauthor{\bsnm{Geyer},~\bfnm{C.}\binits{C.}}
(\byear{1999}).
\btitle{Likelihood inference for spatial point processes}.
In \bbooktitle{Stochastic Geometry ({T}oulouse, 1996)}
(\beditor{O.E. Barndorff-Nielsen}, \beditor{W.S. Kendall}
 and \beditor{M.N.M. van Lieshout},
  eds.).
\bseries{Monogr. Statist. Appl. Probab.}
\bvolume{80}
\bpages{79--140}.
\bpublisher{Chapman \& Hall/CRC, Boca Raton, FL}.
\bid{mr={1673118}}
\bptok{imsref}%
\end{bincollection}
\endbibitem

\bibitem{goulsarkgrab96}
\begin{barticle}[auto:STB|2012/03/20|13:03:14]
\bauthor{\bsnm{Goulard},~\bfnm{M.}\binits{M.}},
  \bauthor{\bsnm{S{{\"a}rkk{\"a}}},~\bfnm{A.}\binits{A.}} \AND
  \bauthor{\bsnm{Grabarnik},~\bfnm{P.}\binits{P.}}
(\byear{1996}).
\btitle{Parameter estimation for marked Gibbs point processes through the
  maximum pseudolikelihood method}.
\bjournal{Scand. J. Stat.}
\bvolume{23}
\bpages{365--379}.
\bptok{imsref}%
\end{barticle}
\endbibitem

\bibitem{JK}
\begin{barticle}[mr]
\bauthor{\bsnm{Jensen},~\bfnm{Jens~Ledet}\binits{J.L.}} \AND
  \bauthor{\bsnm{K{\"u}nsch},~\bfnm{Hans~R.}\binits{H.R.}}
(\byear{1994}).
\btitle{On asymptotic normality of pseudo likelihood estimates for pairwise
  interaction processes}.
\bjournal{Ann. Inst. Statist. Math.}
\bvolume{46}
\bpages{475--486}.
\bid{issn={0020-3157}, mr={1309718}}
\bptok{imsref}%
\end{barticle}
\endbibitem

\bibitem{kuchsore97}
\begin{bbook}[mr]
\bauthor{\bsnm{K{\"u}chler},~\bfnm{Uwe}\binits{U.}} \AND
  \bauthor{\bsnm{S{\o}rensen},~\bfnm{Michael}\binits{M.}}
(\byear{1997}).
\btitle{Exponential Families of Stochastic Processes}.
\bseries{Springer Series in Statistics}.
\baddress{New York}: \bpublisher{Springer}.
\bid{mr={1458891}}
\bptok{imsref}%
\end{bbook}
\endbibitem

\bibitem{lenn24}
\begin{barticle}[auto:STB|2012/03/20|13:03:14]
\bauthor{\bsnm{Lennard-Jones},~\bfnm{J.~E.}\binits{J.E.}}
(\byear{1924}).
\btitle{On the determination of molecular fields}.
\bjournal{Proc. R. Soc. Lond. Ser. A}
\bvolume{106}
\bpages{463--477}.
\bptok{imsref}%
\end{barticle}
\endbibitem


\bibitem{KMM}
\begin{bbook}[mr]
\bauthor{\bsnm{Matthes},~\bfnm{Klaus}\binits{K.}},
  \bauthor{\bsnm{Kerstan},~\bfnm{Johannes}\binits{J.}} \AND
  \bauthor{\bsnm{Mecke},~\bfnm{Joseph}\binits{J.}}
(\byear{1978}).
\btitle{Infinitely Divisible Point Processes}.
\baddress{Chichester}: \bpublisher{Wiley}.
\bid{mr={0517931}}
\bptok{imsref}%
\end{bbook}
\endbibitem


\bibitem{NZ}
\begin{barticle}[mr]
\bauthor{\bsnm{Nguyen},~\bfnm{Xuan-Xanh}\binits{X.X.}} \AND
  \bauthor{\bsnm{Zessin},~\bfnm{Hans}\binits{H.}}
(\byear{1979}).
\btitle{Integral and differential characterizations of the {G}ibbs process}.
\bjournal{Math. Nachr.}
\bvolume{88}
\bpages{105--115}.
\bid{doi={10.1002/mana.19790880109}, issn={0025-584X}, mr={0543396}}
\bptok{imsref}%
\end{barticle}
\endbibitem

\bibitem{ogattane84}
\begin{barticle}[mr]
\bauthor{\bsnm{Ogata},~\bfnm{Yosihiko}\binits{Y.}} \AND
  \bauthor{\bsnm{Tanemura},~\bfnm{Masaharu}\binits{M.}}
(\byear{1984}).
\btitle{Likelihood analysis of spatial point patterns}.
\bjournal{J.~Roy. Statist. Soc. Ser. B}
\bvolume{46}
\bpages{496--518}.
\bid{issn={0035-9246}, mr={0790635}}
\bptok{imsref}%
\end{barticle}
\endbibitem

\bibitem{Preston}
\begin{bbook}[mr]
\bauthor{\bsnm{Preston},~\bfnm{Chris}\binits{C.}}
(\byear{1976}).
\btitle{Random Fields}.
\bseries{Lecture Notes in Mathematics}
\bvolume{534}.
\baddress{Berlin}: \bpublisher{Springer}.
\bid{mr={0448630}}
\bptok{imsref}%
\end{bbook}
\endbibitem

\bibitem{ripl88}
\begin{bbook}[mr]
\bauthor{\bsnm{Ripley},~\bfnm{B.~D.}\binits{B.D.}}
(\byear{1988}).
\btitle{Statistical Inference for Spatial Processes}.
\baddress{Cambridge}: \bpublisher{Cambridge Univ. Press}.
\bid{mr={0971986}}
\bptok{imsref}%
\end{bbook}
\endbibitem


\bibitem{Ruelle70}
\begin{barticle}[mr]
\bauthor{\bsnm{Ruelle},~\bfnm{D.}\binits{D.}}
(\byear{1970}).
\btitle{Superstable interactions in classical statistical mechanics}.
\bjournal{Comm. Math. Phys.}
\bvolume{18}
\bpages{127--159}.
\bid{issn={0010-3616}, mr={0266565}}
\bptok{imsref}%
\end{barticle}
\endbibitem

\bibitem{taka86}
\begin{barticle}[mr]
\bauthor{\bsnm{Takacs},~\bfnm{R.}\binits{R.}}
(\byear{1986}).
\btitle{Estimator for the pair-potential of a {G}ibbsian point process}.
\bjournal{Statistics}
\bvolume{17}
\bpages{429--433}.
\bid{doi={10.1080/02331888608801956}, issn={0233-1888}, mr={0849741}}
\bptok{imsref}%
\end{barticle}
\endbibitem

\end{thebibliography}
\end{document}